%% file: Inspection_games.tex
\documentclass[11pt]{article}
\usepackage{color}
\input{rgb}
\usepackage{amsmath, amssymb,paralist,amsthm}
\usepackage{float,graphicx}
\usepackage{graphics}
\usepackage[format=hang,margin=10pt]{caption}

\def\E{\mathbf{E}}

\def\R{\mathbf{R}}

\def\1{\mathbf{1}}

\def\E{\mathbb{E}}
\def\R{\mathbf{R}}

\newtheorem{prop}{Proposition}[section]
\newtheorem{theorem}{Theorem}[section]

\newtheorem{remark}{Remark}[section]

\numberwithin{equation}{section}

\begin{document}

\centerline{\Large \bf Inspection and crime prevention:}
\smallskip
\centerline{\Large \bf an evolutionary perspective}
\bigskip
\centerline{\bf Vassili Kolokoltsov\footnote{Department of Statistics, University of Warwick, Coventry, CV4 7AL, UK, v.kolokoltsov@warwick.ac.uk}, Hemant Passi\footnote{New Zealand Treasury Wellington, 6011, New Zealand, hemant.passi@treasury.govt.nz}, Wei Yang\footnote{Department of Mathematics and Statistics, University of Strathclyde, Glasgow, G1 1XH, UK, w.yang@strath.ac.uk}}
\bigskip

\begin{abstract}
In this paper, we analyse inspection games with an evolutionary perspective.
In our evolutionary inspection game with a large population, each individual is not a rational payoff maximiser, but  periodically updates his strategy if he perceives that other individuals' strategies are more successful than his own, namely strategies are subject to the evolutionary pressure. We develop this game into a few directions. Firstly,  social norms are  incorporated into the game and we analyse how social norms may influence individuals' propensity to engage in criminal behaviour. Secondly, a forward-looking inspector is considered, namely, the inspector chooses the level of law enforcement whilst taking into account the effect that this choice will have on future crime rates.
Finally, the game is extended to the one with continuous strategy spaces.

\end{abstract}


\section{Introduction}
An inspection game is a non-cooperative game whose players are often called an inspector and an inspectee. It models a situation where the inspectee, which may be an individual, an organisation, a state or a country, is obliged to follow certain regulations but has an incentive to violate them. The inspector tries to minimise the impact of such violations by means of inspections that uncover them. Typically an inspection game has a mixed equilibrium.

Inspection games have a wide variety of applications including arms control \cite{ACKvSZ1996, AC2009, AC2011, H2007}, auditing of accounts \cite{K1968, B1982, B1990}, tax inspection \cite{S1971, R1979, G1984, RW1986, RW1985, AM2004}, environmental protection \cite{R1990, WO1991, GP1992, A1994},  quality control in supply chains \cite{RT1995, TK2007, HL2010}, stockkeeping \cite{FT2008-1} and communication infrastructures \cite{GMCKB2012, CHI2011}. Some of these are surveyed in e.g. \cite{AvS1991, A1986,  ACKvSZ1996, AvSZ2002}. The research on inspection games can contribute to the construction of an effective inspection plan for the inspector when an illegal action is executed strategically.

This subject has been investigated quite extensively during the last five decades. Still, it continues to draw attention. Throughout the 1960s and the three decades that ensured, the underlying motivation was the cold war between the US and the Soviet Union and the desire to monitor the various arms control agreements that were signed by the two superpowers. Analytically, these settings led quite naturally to two-person game formulations with various assumptions about the strategy sets that were feasible to each of the two parties. Since the early 1990s, attention has shifted to other directions and new circumstances require moving from two-person to n-person formulations, from zero-sum games to non-zero-sums.

In the initial work done by Dresher \cite{D1962} and Kuhn \cite{K1963}, they considered sequential inspection games for the treaty of arms reduction where one player wishes to violate the treaty in secret while the other player makes a plan for inspection. In a sequential inspection game, an inspector has to distribute a given number of inspections over a larger number of inspection periods in order to detect an illegal act that an inspectee, who can count the inspector's visits, performs in at most one of these periods. Maschler \cite{M1966} generalized Dresher's modelling.  In Maschler's model, the inspectee (violator) must pay a unit penalty if the violation is exposed by an inspection, but he can escape the exposure by a side payment of penalty $\gamma$. Dresher discusses special cases of $\gamma=1/2$ and $\gamma=1$, and Maschler considered the general case $0\leq \gamma\leq1$. Diamond \cite{D1982} and Ruckle \cite{R1992} also studied two-person zero-sum games between an inspector and an inspectee and a pair of equilibrium strategies was sought.

In the setting of inspections by the customs against smugglers, Thomas and Nisgav \cite{TN1976}  extended Maschler's model to a game in which customs keeps watch on illegal actions of a smuggler by using one or two patrol boats, and the number of boats affects the reward obtained on the capture of the smuggler. They formulated the problem as a multi-stage game, and adopted a numerical method where a one-stage matrix game is repeatedly solved.

Baston and Bostock \cite{BB1991} provided a closed form solution of Thomas and Nisgav's model. They also relaxed the assumption of perfect-capture. The perfect capture assumption means that the inspector certainly captures the inspectee (violators or smugglers) when both players meet. They assumed that the capture probability depends on the number of patrol boats. Gernaev \cite{G1994}  extended Baston and Bostock's work to a model with three patrol boats. Sakaguchi \cite{S1994} and Von Stengel \cite{vS1991} studied inspection games with several smuggling opportunities in a perfect-capture model. Ferguson and Melolidakis \cite{FM1998}  extended further this model by assuming that the smuggler can avoid capture by paying some cost $\gamma \leq1$. If he is captured on his smuggling, the smuggler must pay unit cost.

Since the inspection resources (inspection personnel, equipments or time) are usually limited and complete surveillance is not possible, it is natural to include the first kind error (probability of a false alarm) $\alpha$ given that the inspectee acts legally  and the second kind error (probability of non-detection) $\beta$ that no alarm is raised given that the inspectee acts illegally. This kind inspection game is referred to as an imperfect inspection game, see e.g. Rothenstein and Zamir \cite{RZ2002}, Canty, Rothenstein and Avenhaus \cite{CRA2001}. Rothenstein and Zamir \cite{RZ2002} extended Diamond's  model \cite{D1982} by including fixed errors of the first and second kind for a single inspection.  Canty, Rothenstein and Avenhaus \cite{CRA2001}  solved a non-sequential inspection game with first and second kind errors based on a fixed detection time goal and obtained some special solutions for the sequential case.

It is worth noting that the papers cited above study zero-sum inspection games.
However under certain circumstances, one cannot justify the zero-sum assumption and non-zero-sum formulations are presented, e.g. in Avenhaus and  von Stengel \cite{AvS19912}, Avenhaus, Von Stengel and Zamir \cite{AvSZ2002}, Avenhaus \cite{A2004}, Avenhaus and Canty \cite{AC2005}, Kolokoltsov and Malafeyev \cite{KM2010}, Avenhaus and Kilgour \cite{AK2004}, Hohzaki \cite{H2007} and Deutsch, Golany and Rothblum \cite{DGR2011}.

Avenhaus \cite{A2004}, Kolokoltsov and Malafeyev \cite{KM2010} presented a one-inspector-one-inspectee non-zero-sum one-step imperfect inspection game, where the inspector has a non-detection probability $\beta$. Avenhaus and Canty \cite{AC2005} considered a similar model in a sequential setting, and subject to both kind errors. In these works, they gave an explicit form of the unique Nash equilibrium of this game, i.e. a mixed strategy profile $(p^*, q^*)$, where $p^*$ is the optimal probability with which the inspector inspects and  $q^*$ is the optimal probability with which the inspectee acts illegally.

Avenhaus, Von Stengel and Zamir \cite{AvSZ2002} considered a one-inspector-one -inspectee non-zero-sum imperfect inspection game, where the inspector has a probability of a false alarm $\alpha $ and a probability of non-detection $\beta$, and the inspectee has a behavior strategy $(q, \omega)$, where $q$ is the probability for acting illegally and $\omega$ is a probability distribution on violation procedure given the inspectee acts illegally. They decomposed this game into two games which are solved separately and the solution of one game is a parameter of the second game. First, assuming that the inspectee acts illegally, they formulated a zero-sum game where the inspector faces a hypothesis testing problem $\delta$  and the inspectee chooses an optimal violation mechanism $\omega$. That is, for any test $\delta$, the inspectee will choose a violation mechanism $\omega^*$ to maximise the non-detection probability $\beta(\delta, \omega)$. On the side of the inspector,  first he chooses a fixed false alarm probability $\alpha$ and consider only those tests that result in this error probability. Then he will choose a test $\delta^*$ that minimizes the worst-case non-detection probability $\beta(\delta^*, \omega^*)$. Namely, in equilibrium, the non-detection probability is represented by $\beta(\alpha):=\min_\delta \max_\omega \beta (\delta, \omega)$.
Then they solved the original non-zero-sum game to get the optimal false alarm $\alpha^*$ for the inspector and optimal probability $q^*$ for the inspectee acting illegally.

Avenhaus and Kilgour \cite{AK2004} studied a one-inspector-two-inspectees non-zero-sum inspection game with imperfect inspections. They found a Nash equilibrium where the inspector divides his inspection effort among two inspectees. Further, they showed that when detection functions (probability of detecting), as functions of inspection effort,  are convex, inspection effort should be concentrated on one inspectee chosen at random, but when detection functions are concave it should be spread deterministically over the inspectees.

Hohzaki \cite{H2007}  extended Avenhaus and Kilgour's model \cite{AK2004} to a one-inspector-$n$-inspectees non-zero-sum inspection game, where the inspectee $k$, $k\in[1,n]$, has $l_k$ facilities which are set as targets for the inspection and the inspector has a budget constraint $C>0$. He decomposed this non-zero-sum inspection game into two optimisation problems: an optimal assignment problem of inspection staff and an optimal division problem of the budget. For an optimal assignment problem of inspection staff, he started with assuming resource $y_k$ is assigned to the inspectee $k$, with  $\sum_{k=1}^ny_k=C$. An feasible assignment plan of the inspection staff $x^k=\{x_i^k\geq 0,i=1,..,l_k\}$ must satisfy $\sum_{i=1}^{l_k}c_i^k x_i^k=y_k$, then the inspector aims to minimise the non-detection probability $\beta(y_k)$ over feasible assignment plans $x^k$, $k\in[1,n]$. Further, knowing the minimum non-detection probability $\beta^*(y_k)$, he obtained a Nash equilibrium  $(y^*=\{y_1^*,\dots,y_n^*\}, q^*=\{q_1^*,\dots,q_n^*\})$ of pure strategies of the inspector and mixed strategies of the inspectees by solving an optimal division problem of the budget. Deutsch, Golany and Rothblum \cite{DGR2011} considered a similar non-zero-sum game with single inspector multiple inspectees and  provided closed-form expressions of Nash equilibria.

Gianini, Mayer, Coquil, Kosch and  Brunie \cite{GMCKB2012} generalised the model by Avenhaus \cite{A2004},  Kolokoltsov and Malafeyev  \cite{KM2010} to $m$-inspectors-$n$-inspectees one-step  non-zero sum games. In this game, they assumed the $m$ inspectors are uncoordinated and the $n$ inspectees are non-interacting. Inspectors have probability $p$ to perform the inspection. If an inspector decides for the inspection, then the inspection will be performed on a single randomly chosen inspectee. Given the inspection, each inspectee will have probability $1/2$ to be inspected. The inspectees have respectively probabilities $q_1, \dots, q_n$ of acting illegally. They provided explicit form of mixed strategies $(p^*, q^*_1, \dots, q^*_n)$ for these games.

In the setting of production economics, inspection games are modified in the way that a third player, symbolizing a mediator and coordinating the inspector and the inspectee, is introduced to guarantee that  the inspector and the inspectee behave in accordance with the objectives of the firm. For detailed discussion, see e.g. Borch \cite{B1982}, Anderson and Young \cite{AY1988}, Avenhaus, Von Stengel and Zamir \cite{AvSZ2002}, Fandel/Trockel \cite{FT2008-2}, Fandel/Trockel \cite{FT2011, FT2012}. Another direction of modification of inspection games is to consider the inspector as a Stackelberg leader of the game, see more e.g. in  Andreozzi \cite{A20042, A2008}. McEneaney and Singh \cite{MS2006} studied inspections with spatially nontrivial distributions. A very much related model to the inspection game is the {\it patrolling game}  which involves the scheduling and deployment of patrols, for discussions on the patrolling game see e.g. \cite{AMP2001} and references therein.

Inspection games are also applied to law enforcement.  Tsebelis \cite{T1989} first used inspection games to model phenomena in criminal justice.
According to Tsebelis \cite{T1989}, modifying the size of the penalty does not affect the frequency of crime commitment at equilibrium,
but rather the frequency of law enforcement. Pradiptyo \cite{P2006} refined the inspection game proposed by Tesbelis, by using empirical evidence from various studies, which primarily conducted in the UK, to re-construct the game. They found that an increase in the severity of punishment will reduce the likelihood of the enforcement of the law. Instead of increasing the severity of punishment, theoretically, the authority may reduce individuals' offending behaviour by providing financial compensation to those who do not have criminal records. This is widely known as crime prevention initiative. Interesting psychological analysis of the logic and rationale employed by the players in such game scenarios can be found in Rauhut \cite{R2009}. Altruistic punishment in promoting cooperation in a large population can be found e.g. in \cite{SHTD2011} and \cite{F2005} and many references therein.

In the work by Lindbeck, Nyberg and  Weibull \cite{LNW1999}, they analysed the interplay between social norms and economics incentives in the context of work decisions in the modern welfare state.
 Each individual either works full-time or does not work at all. Each individual chooses to work if and only if his choice results in higher utility than living off the transfer, i.e.
$$u[(1-a)w]>u(b)+\mu-v(x),$$
where $u:[0,\infty)\to [0,\infty)$ is the utility function,  $a\in(0,1)$ is the tax rate, $b\in \R^+$ is the amount of government transfer, $w\in [0,\infty)$ is the wage, $x\in [0,1]$ is the population share living on the transfer, $v(x)$ is a decreasing function describing the disutility from accepting the transfer,  $\mu$ is the utility difference between the leisure of living on the transfer and the intrinsic utility that one may derive from work life. They proved that for every tax rate $a\in(0,1)$, transfer $b\in[0,\infty)$ and expected population share $x\in [0,1]$ of the transfer recipients, there exists a unique critical wage rate $w^*$ such that all individuals with lower wages than $w^*$ choose not to work and those with higher wages than $w^*$ choose to work.  Further, an equilibrium population share $x^*\in [0,1]$ of transfer recipients can be obtained by solving the equation
$$x=\Phi\left (  \frac{u^{-1}[u(b)+\mu-v(x)]}{1-a}\right)$$
where $\Phi:[0,\infty)\to [0,1]$ is a continuously differentiable cumulative probability distribution function of wages for each individual and $u^{-1}$ denotes the inverse function of $u$.
By considering the so-called government balance requirement, that is public transfer spending is equal to public revenues from the wage tax, they showed some interesting observations on the relation between the tax rate $a$ and the transfer $b$.


In this paper, we will investigate a few imperfect inspection games in different settings. First, we show that the canonical two-player inspection game can be generalised to situations where the inspector is responsible for inspecting large populations of interacting individuals whose strategies are subject to an evolutionary pressure. We then examine various extensions, including how social norms  may may influence individuals' propensity to engage in criminal behaviour. Also, we consider how alternative formulations of the inspectorÕs decision problem may affect the aggregate dynamics of crime. Finally, we relax the restriction that individuals face only a binary choice between honesty and crime in order to analysis how our results might generalise to scenarios where individuals can choose the extent of their criminal behaviour.

\section{A conventional inspection game}\label{Canonical games}

In its simplest form, an inspection game is a two-player game between an individual (i.e. inspected) and an inspector. The individual must choose whether to violate or comply with the wishes of the inspector, who himself must decide whether to inspect the individual or not.  Each player makes his strategy choice without observing that of the other. This situation can be illustrated as the following 2x2 normal form game in Table \ref{The canonical two-player inspection game}, where the individual is the row player and the inspector is the column player, and the left (resp. right) entry of each cell corresponds to the payoffs for the individual (resp. inspector).
\begin{table}[ht]
\centering
\begin{tabular}{|l|c|c|c|c|c|}
\hline
 & Inspect & Not Inspect \\[.4em]
\hline
Violate & $-1,1$ & $2,-2$ \\[.4em]
\hline
Comply & $0,-1$ & $0,0$ \\[.4em]
\hline
\end{tabular}
\caption{The simplest two-player inspection game}
\label{The canonical two-player inspection game}
\end{table}

The key feature of an inspection game is that the interests of the individual clash with those of the inspector. The inspector would obviously prefer the individual to comply, but would rather not have to inspect.  However, if the inspector chooses not to inspect, the individual would always prefer to violate. As a result of this conflict, an inspection game contains no pure strategy Nash equilibria: whatever the outcome of the game, either the individual or the inspector would always want to switch their strategy to something that would upset their opponent, who in turn would wish to switch their strategy and so on...

This element of strategic indeterminacy is resolved by resorting to mixed strategies, whereby at least one player randomises over his pure strategies. In fact it is straightforward to show that the unique Nash equilibrium of an inspection game involves both players choosing to play mixed strategies, therefore implying that both inspection and non-compliance occur with positive probabilities.

To describe an inspection game in a general form, we use the following basic notations:
\begin{align*}
r=& \text{legal income received by the individual if he complies}, r>0\\
l=& \text{crimial profit received by the individual if he violates without }\\
&\text{being caught}, l>0\\
f=& \text{fine paid by the individual if he violates and is caught}, f>0\\
c=&\text{cost for the inspector to inspect}, c>0\\
\lambda=& \text{the probability with which the individual is caught if the inspector}\\
&\text{ inspects}, \lambda\in(0,1)
\end{align*}

In a general form of the inspection game between an individual and an inspector can be described as follows: if the individual abides by the law then he receives the legal income $r$, whereas the payoff to the individual that violates depends on the strategy chosen by the inspector. If the inspector did not inspect, the individual escapes with his legal income $r$ plus criminal profits $l$. However if the inspector did inspect, then the individual is either caught with a probability $\lambda$ and has to pay the fine $f$, or escapes with a probability $1-\lambda$ and the criminal profits $l$. The parameter $\lambda$ can also be interpreted as the efficiency of inspection.

Meanwhile, the inspector must incur the fixed cost $c$ in order to inspect. If the individual played Comply, this expense is sunk, and the individual receives no further payoff. However, if the individual chose to violate, then the inspector may catch the individual with the probability $\lambda$ and earn the fine $f$, otherwise the individual escapes, causing the inspector losses $l$. On the other hand, if the inspector chose not to inspect, he sustains losses $l$ when the individual plays Violate, and $0$ when the individual plays Comply.

This game is  summarised by Table \ref{A general inspection game}:
\begin{table}[ht]
\centering
\begin{tabular}{|l|c|c|c|c|c|}
\hline
 & Inspect & Not Inspect \\[.4em]
\hline
Violate & $\big(r+(1-\lambda) l-\lambda f, -c+\lambda f-(1-\lambda) l\big)$ & $(r+l, -l)$ \\[.4em]
\hline
Comply & $(r,-c)$ & $(r,0)$ \\[.4em]
\hline
\end{tabular}
\caption{a general inspection game}
\label{A general inspection game}
\end{table}

In order to ensure that the game can be interpreted as a conventional inspection game, very often the following two additional restrictions on the payoffs are imposed:
\begin{equation}\label{Restriction 1}
  \lambda (l+f)>c
\end{equation}
and
\begin{equation}\label{Restriction 2}
  \lambda f>(1-\lambda)l.
\end{equation}
Restriction \eqref{Restriction 1} states that, given that the individual is violating, the inspector prefers to inspect since the expected benefit from inspecting exceeds the cost. Restriction \eqref{Restriction 2} states that, given the inspector is inspecting, the individual prefers not to violate since the expected fine exceeds the expected criminal profit.

The following theorem gives the unique Nash equilibrium of this game and was  proved in \cite{KM2010} .
\begin{theorem}\label{canonical game equil}
Let $p\in[0,1]$ be the probability with which the individual violates and
$q\in[0,1]$ be the probability with which the inspector inspects.
The unique Nash equilibrium of the two-player inspection game described in Table \ref{A general inspection game} is the mixed strategy profile $(p^*,q^*)$ with
	$$p^*=\frac{c}{\lambda(l+f)}\quad \text{and}\quad q^*=\frac{l}{\lambda(l+f)}.$$
\end{theorem}

\section{An evolutionary inspection game}

In this section, we consider an inspection game with a large population of interacting individuals and a single inspector. In this game, we do not assume that individuals know all relevant features of the game and maximise their payoffs accordingly. Instead, we will replace this conventional payoff maximisation with the idea that individuals modify their strategies after observing the experiences of others. The proportion of individuals who play a particular strategy is then subject to evolutionary pressure over time such that the population share of better performing strategies increases and that of strategies earning lower payoffs decreases. We refer to this model as {\it an evolutionary inspection game}. The key feature of this evolutionary inspection game is that the transmission process by which the strategies of Violate and Comply propagate is fundamentally social in nature.

\subsection{Basic model of an evolutionary inspection game}\label{Evolutionary inspection games}

Formally, in an evolutionary inspection there are $N$ individuals, who choose whether to {\it Violate $(V)$} or {\it Comply $(C)$} and can periodically update their strategy. We assume that the number $N$ is large in order to make a continuous approximation. We call an individual who plays $V$ a {\it violator} and one who plays $C$ a {\it complier}.

The fraction of violators in the population is
$$p:[0,\infty)\to [0,1],\quad t\to p_t.$$
The value $p_t$, called the {\it crime rate} at time $t$, can also be understood as the probability with which each individual violates at time $t$.

The fraction of the population that is inspected by the inspector is
$$q: [0,\infty)\to [0,1] ,\quad t\to q_t.$$
The value $q_t$, called the {\it level of law enforcement} at time $t$, can also be understood as  the probability with which each individual is inspected.

The cost function of law enforcement for the inspector is
$$F:[0,1]\mapsto [0,\infty),\quad q\to F(q),$$
with $F(0)=0$ (that is, the cost of zero inspection is zero).

At any time $t\in[0,\infty)$, a complier always receives the legal income $r$  with certainty, so for any given level of law enforcement $q_t$, a complier expects the payoff:
\begin{equation}\label{payoff of complier-EV}
 \pi_C(t)=:\E [U_I(C,q_t)]=r.
\end{equation}
where  $U_I(C,q)$ denotes the utility of a complier  with any level of law enforcement $q\in[0,1]$.

At any time $t\in[0,\infty)$, if an individual violates, the infraction may go unnoticed with the probability $1-\lambda q_t$ and the violator escapes with the criminal profits $l$; otherwise the violator is caught with the probability $\lambda q_t$ and punished with a fine $f$. Therefore a violator expects the payoff:
\begin{equation}\label{payoff of violator-EV}
 \pi_V(t)=:\E [U_I(V,q_t)]=r+(1-\lambda q_t)l-\lambda q_tf.
\end{equation}
where  $U_I(V,q)$ denotes the utility of a violator  with any law enforcement $q\in[0,1]$.

At any time $t\in[0,\infty)$, roughly speaking there will be $Np_t$ violators, of which a fraction of $\lambda q_t$ is expected to be caught and $1-\lambda q_t$ to escape. Those who are caught have to pay the fine $f$ to the inspector whereas those who escape cause a loss $l$ for the inspector.  Formally, at time $t$, the inspector expects the payoff:
\begin{equation}\label{payoff of inspector-EV}
\E[U_A(q_t,p_t)]=-F(q_t)+Np_t\lambda  q_t f-Np_t(1- \lambda p_t)l
\end{equation}
where $U_A(q,p)$ denoted the utility of the inspector who has a law enforcement level $q\in[0,1]$ and faces a fraction of violators $p\in[0,1]$.

The above analysis can be summarised in Figure \ref{The extensive form representation of the evolutionary inspection game}:

\begin{figure}[ht]
\centering
\resizebox{0.9\textwidth}{!}{\includegraphics{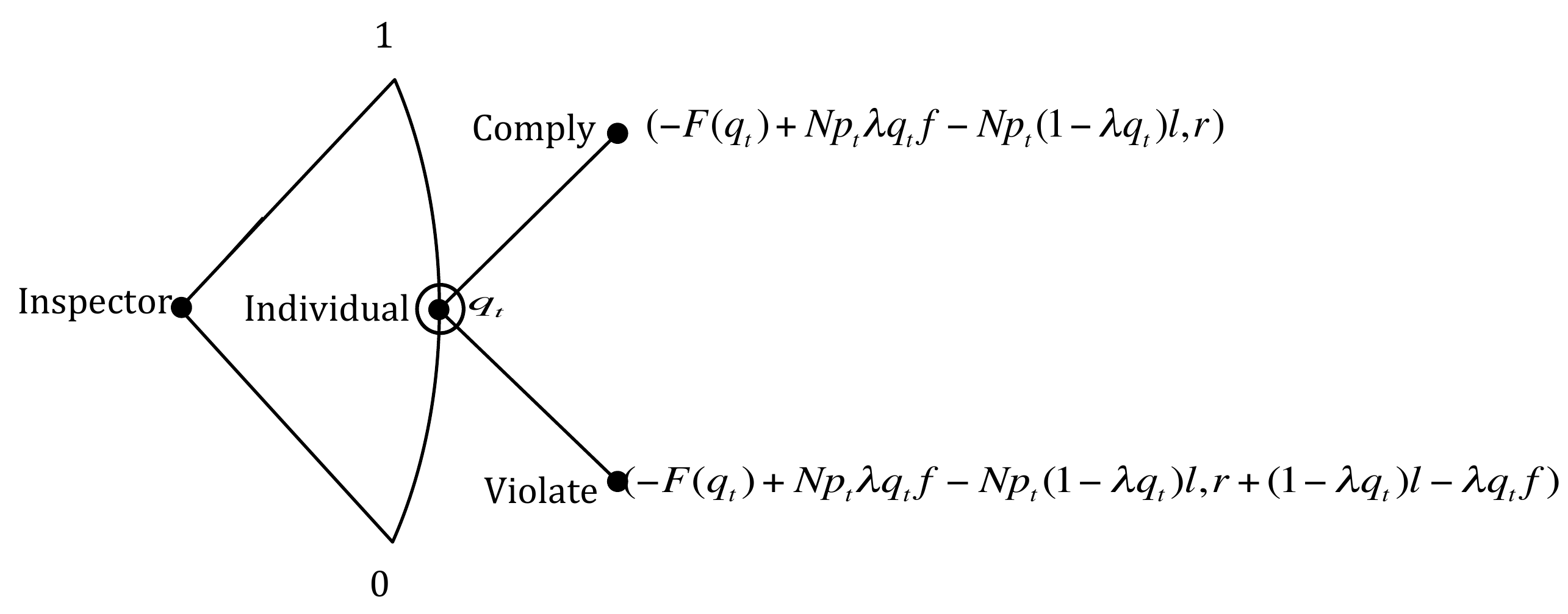} }
\caption{The extensive form of the evolutionary inspection game}
\label{The extensive form representation of the evolutionary inspection game}
\end{figure}

In order to ensure that the game does not degenerate, we impose the  following restriction on the payoff of the individual:
if the inspector inspects with certainty, the individual prefers not to violate, i.e.
{\it 
\begin{equation}\label{Evolutionary-Restriction 1}
\hspace{-10em}\text{Assumption {\bf (A1)}} \quad \quad\quad\quad (1-\lambda )l <\lambda f.\quad\quad\quad\\
\end{equation}}
Inequality \eqref{Evolutionary-Restriction 1} is the same as the restriction \eqref{Restriction 2} in the canonical setting.

In an evolutionary setting, each individual plays $V$ or $C$ as a fixed strategy for a certain period, and then periodically updates his strategy.
At the beginning of each period, some exogenous and fixed fraction $\omega\in [0,1]$ of the population can update their behaviour upon meeting another randomly chosen individual in the population. If two individuals meet and have the same strategy, then that strategy is retained by the updating individual. If however, the two individuals have different strategies, the updating individual may revise his behaviour on the basis of the payoffs enjoyed by the two in the previous period.

More specifically, at a given time $t\in[0,\infty)$, the fraction $1-p_t$ of the population are compliers and $\omega$ of these are in an updating mode and called {\it compliant updaters}. Among these compliant updaters, a fraction $p_t$ of them then meet violators. If the payoff from violating exceeds that from complying, i.e. $\pi_V(t)>\pi _C(t)$, the compliant updater will switch his behaviour with the switching probability
$$\beta(\pi_V(t)-\pi _C(t))$$
where $\beta>0$ is appropriately scaled so that the switching probability  lies in the unit interval $[0,1]$ for any time $t$. If however $\pi_V(t)\leq \pi _C(t)$, the compliant updater does not switch. The coefficient $\beta$ reflects the greater effect on switching of relatively large payoff differences.

The fraction of violators in the next period is equal to the fraction $p_t$ of violators in the previous period, plus the fraction of compliers who convert to violators, and minus the fraction of violators who convert to compliers. Formally the population frequency of violators, at time $t+\Delta t$ ($\Delta t>0$), can be written as
\begin{equation}
 \begin{split}
 \E[  p_{t+\Delta t}]=p_t&+\omega p_t(1-p_t)I_{V>C}\beta(\pi_V(t)-\pi_C(t))\Delta t\\
   &-\omega p_t(1-p_t) I_{C>V}\beta(\pi_C(t)-\pi_V(t))\Delta t
 \end{split}
\end{equation}
where $I_{V>C}$ (resp. $I_{C>V}$) is an indicator function that equals to one if $\pi_V(t)>\pi_C(t)$ (resp. $\pi_C(t)>\pi_V(t)$) and zero otherwise. Since the population is sufficiently large, we can replace $\E[p_{t+\Delta t}]$ by $p_{t+\Delta t}$ and get
\begin{equation}\label{p(t+dt)}
   p_{t+\Delta t}=p_t+\omega\beta p_t(\pi_V(t)-\bar \pi(t))\Delta t\\
\end{equation}
where
\begin{equation}\label{bar pi(t)}
\bar\pi (t)=p_t\pi_V(t)+(1-p_t)\pi_C(t)
\end{equation}
 is the average payoff at time t for the population as a whole. Subtracting $p_t$ from both sides of equation \eqref{p(t+dt)}, dividing by $\Delta t$ and taking the limit as $\Delta t\to 0$, one gets the well-known {\it replicator equation}  (cf. e.g. Taylor and Jonker \cite{TJ1987}  or Zeeman \cite{Z1979})
\begin{equation}\label{replicator eq}
  \dot p_t= \frac{d p_t }{dt}= \omega \beta p_t(\pi_V(t)-\bar \pi(t)).
\end{equation}
This equation states that the growth rate of the fraction of violators is proportional to the amount by which the payoff from violating exceeds that from complying.

By \eqref{payoff of violator-EV}, \eqref{bar pi(t)} and \eqref{payoff of complier-EV}, we can rewrite \eqref{replicator eq}  and derive a replicator equation for the frequency of violators in the population, which is presented in the following.
\begin{prop}
In an evolutionary inspection game with a large population of individuals, the replicator equation for the frequency of violators is given by
  \begin{equation}\label{replicator eq spe}
   \dot p_t= \omega \beta p_t(1-p_t) (l-\lambda  q_t(l+f)).
  \end{equation}
\end{prop}

In this evolutionary inspection game, the replicator dynamic \eqref{replicator eq spe} is not a best-response dynamic. That is, at any time $t$, individuals do not adopt the strategy that maximises their payoffs with regard to the prevailing level of crime $p_t$ and law-enforcement $q_t$. In fact, individuals need not even be aware of the overall levels of crime and law-enforcement. Instead, individuals periodically meet with other members of the population, observe their strategies, and consider updating their behaviour if they perceive that other individuals' strategies are more successful than their own.

The implications of the replicator dynamic are that although no individual can be considered to be a rational payoff maximiser, strategies are subject to the evolutionary pressure. In other words, strategies that are successful in terms of yielding above-average payoffs will propagate through the population by imitation. Therefore the replicator equation is said to be a payoff-monotonic dynamic.

Next, we will analyse how the inspector responds dynamically to the evolving level of crime $p_t, t\geq 0$. In this basic model of an evolutionary inspection game, we assume that the crime rates $p_t$ are visible for the inspector and that the inspector plays a straightforward best-response to the visible crime rates. More specifically, at each time $t\geq 0$, the inspector takes the crime-rate $p_t$ as given and tries to find a level of law enforcement $q_t$ so as to maximise his expected payoff. Specifically, at each time $t\geq 0$ and given a crime rate  $p_t$, the inspector aims to maximise his payoff \eqref{payoff of inspector-EV}, i.e.
\begin{equation}\label{objective function of inspector}
  \max_{q_t\in[0,1]}\left\{-F(q_t)+Np_t\lambda  q_tf-Np_t(1-\lambda  q_t)l\right\}.
\end{equation}

\begin{remark}
One may think that this straightforward best-response strategy is somewhat naive as it assumes that the inspector does not choose $p_t$ strategically so as to influence future levels of crime rate. However, one can also consider the situation that, perhaps due to political or resource issues,  the inspector may not trade-off non-best response strategy in the short-run in the pursuit of some long-run objectives. At any case, this simplifying assumption provides an interesting benchmark from which other assumptions can be compared, see e.g. subsection \ref{Forward-looking inspectors}.
\end{remark}

{\it Assumption {\bf (A2)}:  The cost function $F$ is twice continuously differentiable and
$$0<F'(q)\leq N\lambda  (l+f) \text{ and } F''(q)>0,\quad \text{for all } q\in[0,1],$$
where $F'$ and $F''$ denote the first-order and the second-order derivatives of $F$, respectively.}

By assumption {\bf (A2)}, for each $t\geq 0$, the maximum of \eqref{objective function of inspector} is achieved only at one point. Let the unique best-response for the inspector to any crime rate $p\in[0,1]$ is denoted by
\begin{equation}\label{best respond}
  \hat q [p]:=\arg\max_{q\in[0,1]}\{-F(q)+Np\lambda  qf-Np (1-\lambda  q)l\}.
\end{equation}
Define
$$\bar p:=\min\{p\in[0,1]:\hat q [p]=1\}.$$
Clearly for $p\in[\bar p, 1]$,  we have $\hat q (p)=1$. Essentially, the assumption {\bf (A2)} requires that there exists a critical level of crime rate $\bar p\in(0,1]$ such that if the prevailing crime rate is higher than $\bar p$, the inspector will inspect with certainty.

A very basic example of the cost function $F$ in our setting is
\begin{equation}\label{example F}
F(q)=\alpha q^2
\end{equation}
with a constant $0<\alpha \leq \frac{N\lambda  (l+f)}{2}$. In this case, the unique best-response for the inspector to crime rate $p\in[0,1]$ is given by
$$\hat q [p]=\min\left\{\frac{N\lambda (f+l)}{2 \alpha}p, 1\right\}$$
which is a linear function in $p$. The critical level of crime rate $\bar p$ is given by
$$\bar p= \frac{2 \alpha}{N\lambda (f+l)}.$$
We will keep this example in mind in the following analysis.

\begin{prop}\label{q monoto in p}
Let assumptions {\bf (A1)} and {\bf (A2)} hold. Then  the best-response function for the inspector
$$\hat q: (0,\bar p)\mapsto (0,1), \quad p\mapsto \hat q(p),$$
defined in \eqref{best respond}, is strictly increasing in $p$.
\end{prop}

\proof

By the definition of the best response function $\hat q(\cdot)$ in \eqref{best respond}, we have  the following first order condition
\begin{equation}\label{first order condition}
  -F'(\hat q(p))+Np\lambda  (l+f)=0
\end{equation}
and the second order condition
\begin{equation}\label{second order condition}
  -F''(\hat q(p))< 0
\end{equation}
are satisfied. Differentiating \eqref{first order condition} with respect to $p$, together with \eqref{second order condition}, yields
$$\frac{d\hat q}{dp}(p)=\frac{N\lambda  (l+f)}{F''\left(\hat q(p)\right)}> 0$$
which completes the proof.
\qed

Proposition (\ref{q monoto in p}) says that the inspector's optimal level of law enforcement is increasing in the crime rate.

After having discussed the dynamics for both individuals and the inspector, we are now in a position to analyse the qualitative behaviour of the system. Substituting \eqref{best respond} into \eqref{replicator eq spe}, one gets that the fixed points of the replicator equation \eqref{replicator eq spe} occur at  $p^*=0$,  $p^*=1$,  and $p^*$ such that
\begin{equation}\label{q optimal}
\hat q (p^*)=\frac{l}{\lambda (l+f)}.
\end{equation}
Clearly, by the assumption {\bf (A1)},   $\hat q (p^*)$ in \eqref{q optimal} is within $(0,1)$. Since we have proved that $\hat q(p)$ is strictly increasing in $p$ in Proposition \ref{q monoto in p} , the solution $p^*\in (0,1)$ of \eqref{q optimal} is unique, we refer to this fixed point $p^*\in(0,1)$ as the {\it interior fixed point} of \eqref{replicator eq spe}.

\begin{theorem}\label{evolutionary equlibrium}
Let assumptions {\bf (A1)} and {\bf (A2)} hold. Then
 \begin{enumerate}
\item[{\rm (i)}] The fixed points $p^*=0$ and $p^*=1$ of  \eqref{replicator eq spe} are unstable.
\item[{\rm (ii)}] There exists an unique stable interior equilibrium of the evolutionary inspection game $(p^*,q^*)$ with
 $$p^*=\frac{1}{N\lambda (l+f)} F'(\frac{l}{\lambda (l+f)})\quad \text{ and } \quad q^*=\frac{l}{\lambda (l+f)}.$$
  \end{enumerate}
\end{theorem}

\proof
To analyse the stability of the fixed points of this system, we apply the method of linearizing around these  fixed points (c.f. \cite{S1994}). Specifically,
define a function  $j:[0,1]\mapsto \R$ by
\begin{equation}\label{j(p)}
  j(p)=\omega\beta p(1-p) \big(l-\lambda \hat q (p)(l+f)).
\end{equation}
Let $p^*$  denote the fixed point of the system and $h_t:=p_t-p^*$ be a small perturbation away from $p^*$. The linear equation
$$\dot h_t=h_t\frac{d j}{dp}(p^*)$$
is called the linearization about the point $p^*$. The important fact is that the sign of $\frac{dj}{dp}(p^*)$ tells us the stability of the fixed point $p^*$, i.e., $p^*$ is stable if $\frac{dj}{dp}(p^*)<0$ and unstable if $\frac{dj}{dp}(p^*)>0$. Differentiating the function $j(p)$ \eqref{j(p)} with respect to $p$ gives us
$$\frac{dj}{dp}(p)=\omega\beta (1-2p)\left (l-\lambda \hat q(p)(l+f)\right)-\omega\beta \lambda \frac{d\hat q(p)}{dp}p(1-p)(l+f).$$

 At the point $p^*=0$, we obtain
$$\frac{dj}{dp}(0)=\omega\beta [l-\lambda \hat q(0)(l+f)]=l>0.$$
Therefore the fixed point $p^*=0$ is unstable.

In order to check whether there exists an interior point $p^*\in(0,1)$ and to determine the its' value, we substitute \eqref{q optimal}  into the authorities' first order condition \eqref{first order condition}. Together with the assumption {\bf (A2)} we get
 $p^*=\frac{F'(q^*)}{N\lambda (l+f)}\in(0,1)$ where $q^*= \frac{l}{\lambda (l+f)}\in(0,1)$ . For the stability of this interior fixed point, we have
$$\frac{dj}{dp}(p^*)=-\omega\beta  \lambda \frac{d\hat q(p^*)}{dp } p^*(1-p^*)(l+f)<0$$
since $\hat q$ is increasing in $p$, i.e. $\frac{d\hat q(p)}{dp }>0$ for any $p\in(0,\bar p)$. Therefore, this interior fixed point is stable.

 At the fixed point $p^*=1$, we have
$$\frac{dj}{dp}(1)=-\omega\beta  [l-\lambda \hat q(1)(l+f)].$$
Since there exists an interior fixed point $p^*\in(0,1)$ such that $\hat q (p^*)=\frac{l}{\lambda (l+f)}\in(0,1)$, we have $\hat q (1)\in(\frac{l}{\lambda (l+f)},1]$, by Proposition \ref{q monoto in p}.  Therefore $\frac{dj}{dp}(1)>0$, which implies the fixed point $p^*=1$ is unstable.

\qed

\begin{remark}
The mixed strategy equilibrium $(p^*,q^*)$ in Theorem \ref{evolutionary equlibrium} is essentially the same as the one derived in Theorem \ref{canonical game equil}, in the canonical setting. The main difference can be traced to the fact that the role of the authorities' strategy variable  $q$ differs in the two games. In Section \ref{Canonical games}, the authorities faced a binary choice between inspecting and not inspecting, and $q$ simply specified the authorities' randomisation between these two choices. Given such a randomisation, the authorities' expected cost of inspection would be $cq$, which amounts to a linear cost function with constant marginal costs $c$. In the evolutionary model, however, $q$ is interpreted as a continuous choice variable with a more general cost function $F(q)$ for which equilibrium marginal costs are $F'(q^*)$. Comparing the expressions for $p^*$ and $q^*$ in Theorem \ref{canonical game equil} and Theorem \ref{evolutionary equlibrium}, we can therefore see that their interpretations are essentially the same.
\end{remark}

\subsection{Crime and social norms}

The key feature of an evolutionary inspection game is that the process of individual behavioural change is fundamentally social in nature. Given this position, it is not difficult to imagine applications where other social factors may also exert an influence on individuals' decisions.

Abiding by the law is a social norm but the strength of this norm depends on the size of the population share that adheres to the law. A very low overall crime rate helps to foster a culture of honesty in the population which serves to reinforce individuals' sense of moral responsibility. Conversely if crime is rife,  individuals may begin to feel that breaking the law is socially acceptable.

In this subsection, we will incorporate social norms into our model and analyse how social norms may influence individuals' propensity to engage in criminal behaviour. Inspired by the work \cite{LNW1999}, we model social influence as the degree of disutility (or discomfort) experienced by an individual from choosing a dishonest action. The larger the share of the population that is law-abiding, the more intense is the social discomfort from crime. Therefore the intensity of the social norm is determined endogenously in the model and the degree of disutility is a function of the share of the population $p$. Formally, the disutility function is
$$g: [0,1]\to [0,\infty), \quad p\to g(p).$$

The model described in the subsection \ref{Evolutionary inspection games} should be modified in the way that the disutility $g(p_t)$ is subtracted from \eqref{payoff of violator-EV}.
Then the replicator equation with social norm is given by
\begin{equation}\label{p+social norm}
\dot p_t = \omega \beta p_t(1-p_t) \big(l-\lambda q_t(l+f)-g(p_t)\big).
\end{equation}

As in  subsection \ref{Evolutionary inspection games}, we have the fixed points of \eqref{p+social norm} at the boundaries $p^*=0$ and $p^*=1$. However, the fixed point equation for interior equilibria $p^*\in(0,1)$ now is such that
\begin{equation}\label{social norm equilibrium}
  l-\lambda  \hat q (p^*)(l+f)=g(p^*)
\end{equation}
where $\hat q (\cdot)$ is given by \eqref{best respond}. Define a function $W: [0,1]\to \R$ by
\begin{equation}\label{W1}
W(p)= l-\lambda \hat q (p)(l+f).
\end{equation}

\begin{figure}[ht]
\centering
\resizebox{0.8\textwidth}{!}{\includegraphics{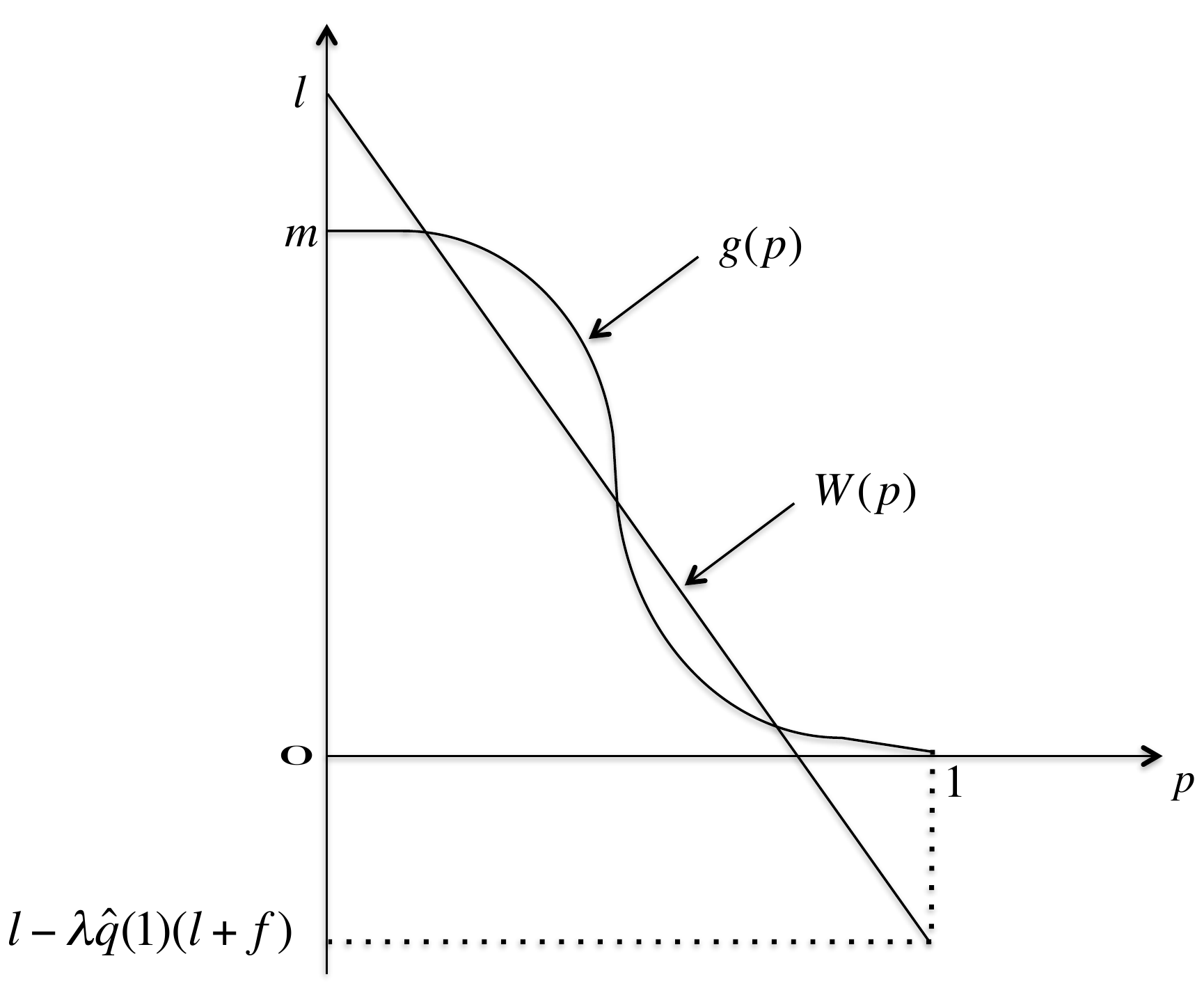} }
\caption{The interior fixed points in the presence of social norms}
\label{The fixed point equation in the presence of social norms}
\end{figure}

{\it Assumption {\bf (A3)} The disutility function $g(p)$ is monotonically decreasing in $p$. It decreases rapidly from a high to a low value at some intermediate value of $p$, namely it takes a sigmoid shape (as illustrated in Figure \ref{The fixed point equation in the presence of social norms}). Let $g(0)=m$, with $0<m<l$, and $g(1)=0$.
}

The sigmoid shape of the disutility function $g(\cdot)$ captures the idea that aggregate behaviour can encourage both honest and dishonest behaviours at the individual level.

\begin{prop}\label{SN}
 Let assumptions {\bf (A1)} and {\bf (A3)} hold. Then the fixed points $p^*=0$ and $p^*=1$ of \eqref{p+social norm} are unstable.
\end{prop}
As a direct consequence of Proposition \ref{SN}, there exists at least one stable interior fixed point of \eqref{p+social norm}.

\proof
Following the same procedure as in the proof of Theorem \ref{evolutionary equlibrium}, the dynamics of the system are now governed by:
\begin{equation}
   \dot p_t = \tilde h(p_t)=\omega \beta p_t(1-p_t) \big(l-\lambda \hat q(p_t)(l+f)-g(p_t)\big).
  \end{equation}
Differentiating the function $\tilde h$ with respect to p
$$\frac{d \tilde h}{d p}(p)=\omega \beta (1-2p)[l-\lambda \hat q(p)(l+f)-g(p)]-\omega \beta p(1-p)\left(\lambda \frac{d\hat q(p)}{dp}(l+f)+g'(p)\right).$$
Then we evaluate this derivative at the boundary fixed point $p^*=0$ and obtain:
$$\frac{d \tilde h}{d p}(0)=\omega \beta(l-\theta)>0.$$
Therefore, $p^*=0$ is an unstable fixed point.

 The left-hand-side of \eqref{social norm equilibrium} for interior fixed points is the function $W$ defined in \eqref{W1}. The function $W$ is decreasing in $p$, as it is simply a linear transformation of the authorities' best-response function $\hat q(\cdot)$ which is increasing in $p$ by Proposition \ref{q monoto in p}.
Consequently, the non-linearity of the disutility function $g$ allows for the possibility of multiple interior equilibria. Rewriting \eqref{social norm equilibrium} yields
$$\hat q(p^*)=\frac{l-g(p^*)}{\lambda (l+f)}=1-\frac{g(p^*)+\lambda f-(1-\lambda)l}{\lambda(l+f)}.$$
Since $g(0)=m<l$, $g(1)=0$ and $g(p)$ is decreasing, we have $0<g(p^*)<l$. Together with the restriction \eqref{Evolutionary-Restriction 1}, we have
$$0<\frac{g(p^*)+\lambda f-(1-\lambda)l}{\lambda(l+f)}<1,$$
namely, there exists an   interior fixed point $p^*\in(0,1)$ of \eqref{p+social norm} such that
$$0<\hat q(p^*)<\frac{l}{\lambda (l+f)}.$$
Therefore we have, at the boundary fixed point $p^*=1$
$$\frac{d \tilde h}{d p}(1)=-\omega \beta (l-\lambda \hat q(1)(l+f))>0$$
that is $p^*=1$ is an unstable fixed point.
\qed

\begin{remark}
  The inclusion of social norms provides a mechanism for positive feedback. When the initial rate of crime is low, the social disutility from crime is high (i.e. the social norm is strong), which encourages more individuals to behave honestly. This in turn acts to further reduce the rate of crime and strengthen the norm, and so on. Conversely, when the rate of crime is high, the social norm is weak and therefore the process acts in reverse.
\end{remark}

Significantly, the presence of multiple equilibria allows us to conclude that the inclusion of social norms implies that societies with very similar social norm and level of punishment can experience very different equilibrium crime-rates if their histories are different.

To illustrate this point, we will take a specific example of the cost function $F$ and analysis how changes in the fine parameter $f$ might affect the qualitative properties of these equilibria. Let $F$ be of the form of \eqref{example F} with $\alpha =\frac{1}{2}N\lambda (l+f)$, i.e.
$$F(q)=\frac{1}{2}N\lambda (l+f)q^2.$$
Then the best respond $\hat q$ defined in \eqref{best respond} is
$$\hat p(p)=p, \quad \text { for all } p\in[0,1]$$
and in turn the function $W$ in\eqref{W1} is
\begin{equation}\label{W2}
W(p)=l-p\lambda (l+f).
\end{equation}
When $p=1$, $W(1)=l-\lambda (l+f)<0$ by the assumption {\bf (A1)}. In this example, the solutions of \eqref{social norm equilibrium} are the intersection points of the linear function $W(p)$ in \eqref{W2} with the nonlinear disutility function $g(p)$ with $g(0)=m<l$ and $g(1)=0$.

\begin{figure}[ht]
\centering
  \includegraphics[width=1\linewidth]{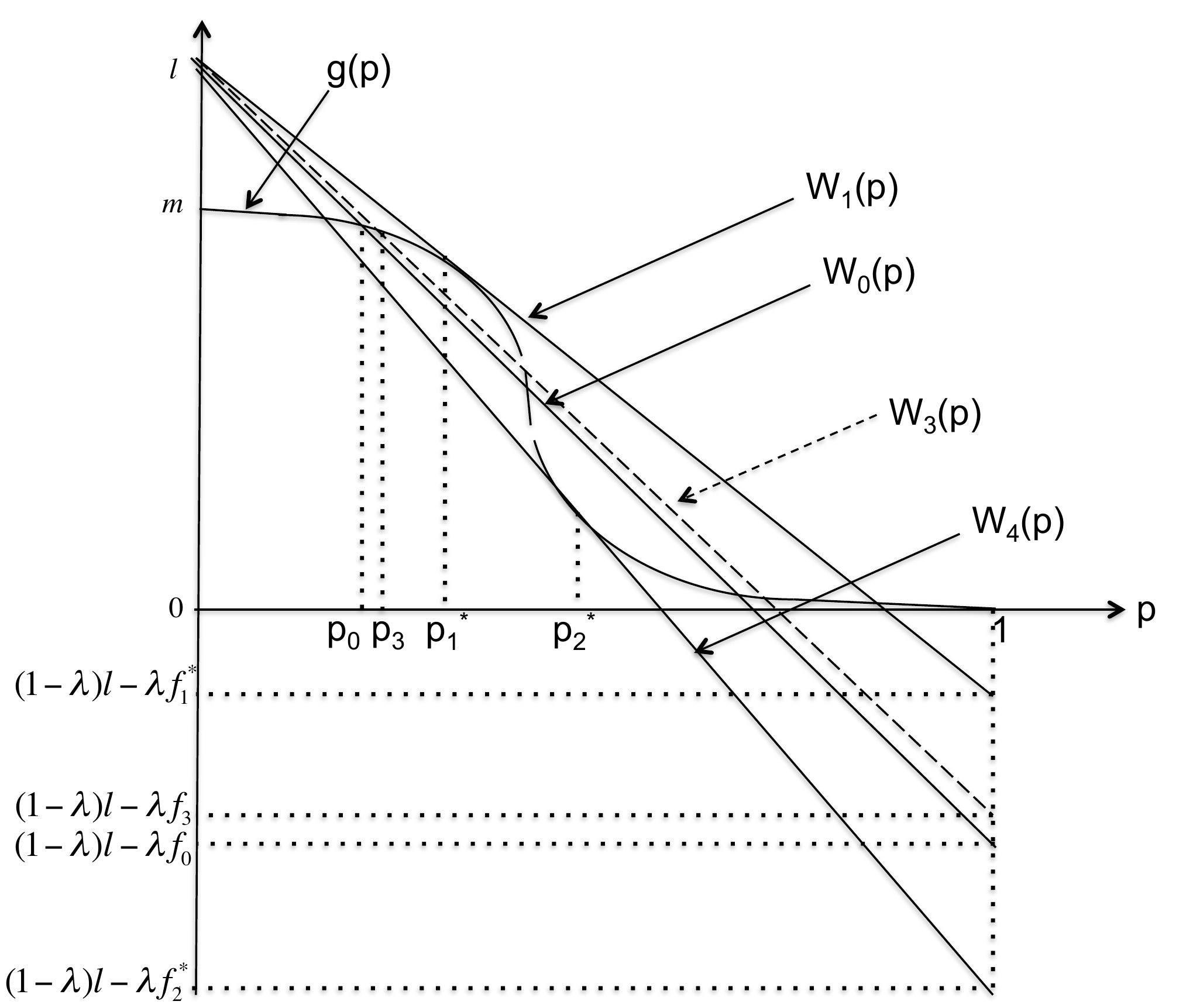}
  \caption{Fine and equilibria}
\label{Fine}
\end{figure}

Let's illustrate the effect of fine changes on the equilibrium with the help of Figure \ref{Fine}. We start with the situation that a society is initially at a low-crime equilibrium $p_0$ when the fine is $f_0$ (corresponding to the line $W_0$). From this point, successive policy changes can allow the inspector to take a more lenient view on crime, namely to decrease the fine from $f_0$ to $f_3$ (i.e. move line $W_0$ to line $W_3$) without experiencing much increase in the equilibrium crime rate, i.e. the equilibrium crime rate increases slightly from $p_0 $ to $p_3$. However, if the level of fine falls below the critical value  $f_1^*$ (i.e. the line $W_0$ is moved above the line $W_1$), there will be a drastic increase in the rate of crime (i.e. moving from $p^*_1$ towards $p^*_2$) as the honest norm is 'destroyed'. Moreover, undoing the policy change by increasing fine from lower than $f^*_1$ to $f_0$ (i.e. by moving a line above the line $W_1$ back to $W_0$) will not be sufficient to restore the low crime equilibrium $p_0$, since by being excessively lenient, the inspector has unwittingly fostered an undesirable social norm. Therefore the level of fine must now be increased beyond $f_2^*$ in order to restore the population to honesty and to restore a low crime rate.

\begin{figure}[ht]
\centering
  \centering
  \includegraphics[width=.7\linewidth]{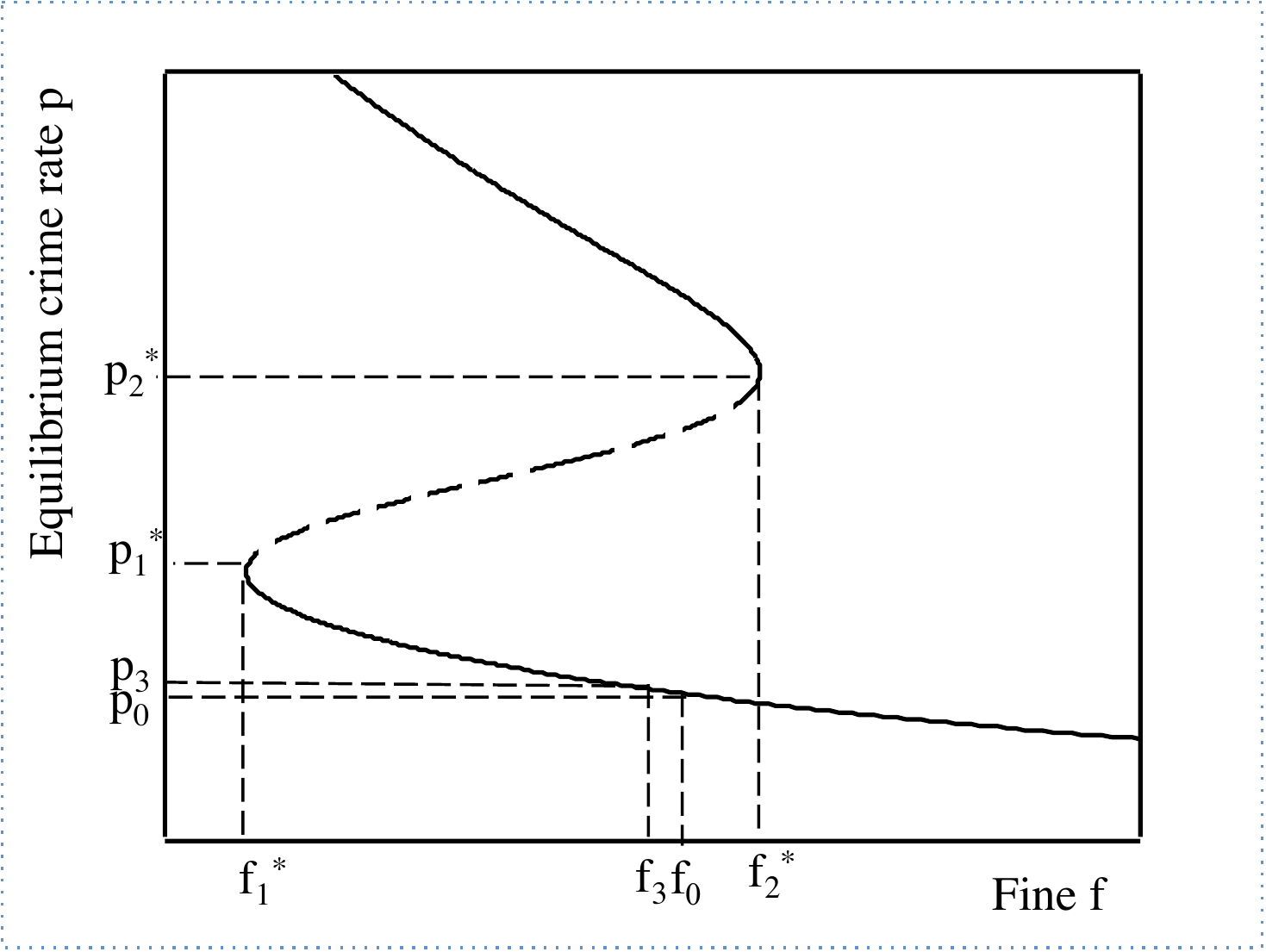}
  \caption{Bifurcation diagram}
\label{Bifurcation diagram}
\end{figure}
The bifurcation curve in Figure \ref{Bifurcation diagram} summaries  the relation between the fine parameter $f$ and the equilibrium rate of crime $p$.
As a result, we conclude that in the presence of social norms, seemingly minor policy changes may cause discontinuous changes to the overall level of crime. Moreover, when such changes in the crime rate occur, they are likely to be highly path-dependent.

\begin{remark}

The multiple projection of the equilibrium curve on the axis of the basic parameter of the model is a general phenomenon in real life. For example, in economics Bowles \cite{B2004} examined market's instabilities of equilibrium prices. In biology, one observes the situations when small changes in drug consumptions can lead to basically irreversible changes in equilibrium states of an organism, see e.g. \cite{GPST2011, GST2010}.
\end{remark}

\subsection{Forward-looking inspectors}
\label{Forward-looking inspectors}

Up until now, we have assumed that the inspector observes the crime rate and simply plays his best-response level of law-enforcement at each time $t$. In this section we revise this assumption by considering the optimal strategy for an inspector who is more {\it forward-looking}. Specifically, we consider the problem faced by an inspector who would like to choose the level of law-enforcement whilst taking into account the effect that this choice will have on future crime rates.

This can be analysed formally by replacing the static maximisation problem for the inspector as described in  \eqref{objective function of inspector} with an optimal control problem. Rather than optimising at each time $t$, the inspector now chooses an entire trajectory of the law enforcement variable $\{q.\}=\{q_t, t\geq 0\}$ so as to 'control' the crime rates $p_t$ via the replicator equation \eqref{replicator eq spe}. In the terminology of optimal control problems, $q_t$ is referred to as the  {\it control variable} and $p_t$ as the {\it state variable}.  The objective of the inspector is to choose the trajectory $\{q.\}$ (i.e. a policy of law enforcement) that maximises the total discounted payoffs with an infinite planning horizon and a discount rate $\delta>0$. In this subsection, social norms do not play a role in the game.

Formally, the inspector solves the following optimal control problem
$$\max_{\{q.\}} \int^{\infty}_0 e^{-\delta t}\big(-F(q_t)+Np_t\lambda q_tf-Np_t(1-\lambda q_t)l\big)dt$$
subject to
\begin{equation}\label{p evolution}
\dot p_t=\omega \beta p_t(1-p_t)(l-\lambda q_t(l+f)):=P(p_t,q_t)
\end{equation}
with an initial value $p_0\in[0,1]$. The Hamiltonian function of this optimal control problem $H:[0,1]\times[0,1]\times [0,\infty)\times\R\to \R$ is written as:
$$H(p,q,t,\mu)=e^{-\delta t}\left(-F(q)+Np\lambda qf-Np(1-\lambda q)l\right)+\mu \omega \beta p(1-p)(l-\lambda q(l+f))$$
where $\mu$ is the dynamic Lagrange multiplier.
By the {\it Maximum Principle}, the first order conditions are:
\begin{align}\label{H_q}
&\frac{\partial H}{\partial q}(p_t,q_t,t,\mu_t)=\notag\\[.4em]
&e^{-\delta t}\left(-F'(q_t)+Np_t\lambda (l+f)\right)-\mu_t \omega \beta p_t (1-p_t)\lambda (l+f)=0
\end{align}
and
\begin{align}\label{H_p}
&\frac{\partial H}{\partial p}(p_t,q_t,t,\mu_t)=\notag\\[.4em]
&e^{-\delta t}N\left(\lambda q_tf-(1-\lambda q_t)l\right)+\mu_t\omega \beta (1-2p_t)(l-\lambda q_t(l+f))=-\dot \mu_t.
\end{align}
Next, we will manipulate  \eqref{H_q} and \eqref{H_p} to eliminate $\mu$ and $\dot \mu$ and get an evolution equation for $q_t$. Differentiating \eqref{H_q} with respect to $t$ yields
\begin{equation}\label{Dt}
\begin{split}
-\delta &e^{-\delta t} \big(-F'(q_t)+Np_t\lambda (l+f)\big)+ e^{-\delta t} \left(-F''(q_t)\dot q_t+N\dot p_t\lambda (l+f)\right)\\
&-\dot \mu_t \omega \beta p_t (1-p_t) \lambda (l+f)- \mu_t\omega \beta (1-2 p_t)\dot p_t \lambda (l+f)=0.
\end{split}
\end{equation}
Then we solve for $\dot \mu_t$ from \eqref{H_p} and plug it into \eqref{Dt}. Together with the dynamics $\dot p_t $ in \eqref{p evolution}, the terms which contain $\mu_t$ are cancelled. Then we get
$$-\delta \big(-F'(q_t)+Np_t\lambda (l+f)\big)-F''(q_t)\dot q_t=0.$$
Hence, we derive the following equation for the evolution of $q_t$
\begin{equation}\label{q evolution}
\dot q_t=\frac{\delta\big(F'(q_t)-Np_t\lambda (l+f)\big)}{F''(q_t)}:= Q(p_t,q_t).
\end{equation}
Equations \eqref{p evolution} and \eqref{q evolution} constitute a pair of first-order differential equations, which describe the optimal trajectories for $p_t$ and $q_t$.

\begin{theorem}
\begin{enumerate}
Let assumptions {\bf (A1)} and {\bf (A2)} hold. Then,
\item[{\rm (i)}]
 there exists only one interior fixed point $(p^*,q^*)$ of the equations \eqref{p evolution} and \eqref{q evolution} such that
\begin{equation}\label{L1}
(p^*,q^*)=\left(\frac{F'(\frac{l}{\lambda (l+f)})}{N\lambda (l+f)}, \frac{l}{\lambda (l+f)}\right)
\end{equation}
\item[{\rm (ii)}] the fixed point $(p^*,q^*)$ is an unstable saddle point.
\end{enumerate}
\end{theorem}

\proof

 The fixed points of the equations \eqref{p evolution} and \eqref{q evolution}  can be found by  examine the {\it nullclines}, that is the sets of points at which both $\dot p_t$ and $\dot q_t$ equal zero. In other words, any fixed point $(p^*,q^*)$ satisfies that
$$P[p^*,q^*]=Q[p^*,q^*]=0.$$
Specifically, let $P[p^*,q^*]=0$ together with assumption {\bf (A1)}, we get
$$ p^*=0, \quad p^*=1 \quad \text{and } p^*\in(0,1) \, \text{such that } \, q^*[p^*]=\frac{l}{\lambda (l+f)}\in(0,1).$$
Let $Q(p^*,q^*)=0$ we get
$$p^*=\frac{F'(q^*)}{N\lambda (l+f)}\in(0,1) \text{ with } q^*\in(0,1)$$
since by assumption {\bf (A2)}, for any $p\in(0,1)$, $0<F'(\frac{l}{\lambda (l+f)})<N\lambda (l+f)$.
Therefore, there exists an unique fixed point of the system equations \eqref{p evolution} and \eqref{q evolution}, and the fixed point is at
$$(p^*,q^*)=\left(\frac{F'(\frac{l}{\lambda(l+f)})}{N\lambda (l+f)}, \frac{l}{\lambda (l+f)}\right).$$

(ii) The stability of this fixed point can be analysed by examine the signs of the eigenvalues of the linearised system of differential equations \eqref{p evolution} and \eqref{q evolution} (c.f. \cite{S1994}).

Linearising the system about the fixed point $(p^*,q^*)$ \eqref{L1} gives:
\[
\begin{bmatrix}
 \dot p_t\\
  \dot q_t
 \end{bmatrix}
 =
 A
  \begin{bmatrix}
  p_t-p^*\\
 q_t-q^*
 \end{bmatrix}
 \]
 where $A$ is a $2\times 2$ Jacobian matrix at the fixed point $(p^*,q^*)$, i.e.
 \[
 A=\begin{bmatrix}
  \frac{\partial P}{\partial p}&\frac{\partial P}{\partial q}\\[.4em]
  \frac{\partial Q}{\partial p}& \frac{\partial Q}{\partial q}
 \end{bmatrix}_{(p^*,q^*)}
 =
 \begin{bmatrix}
  0&-\omega \beta p^*(1-p^*)\lambda (l+f)\\
  -\frac{\delta N\lambda (l+f)}{F''(q^*)}& \delta.
 \end{bmatrix}
 \]
The determinant of $A$ is
$$det(A)=-\omega \beta p^*(1-p^*)\lambda (l+f) \frac{\delta N\lambda (l+f)}{F''(q^*)}<0$$
since the cost function $F$ is strictly convex, by the assumption {\bf (A2)}. Hence the fixed point $(p^*,q^*)$ is a saddle point. Since the trace of $A$ is $Tr(A)= \delta>0$, this point is unstable.
\qed


\begin{remark}
In this work we concentrate on the case with strictly convex cost function $F(q)$.  In the case where the equilibrium costs are locally concave (i.e. $F''(q^*)<0$), it is possible to show that the eigenvalues of the system are complex, with a positive real component. Therefore the interior steady-state is a spiral node; all nearby trajectories spiral away  from the steady-state.   In this case, the precise optimal control is difficult to characterise in general terms, and depends on whether there are any additional steady-states at the boundaries which, providing they exist, will always be saddle points. However at the very least, we can conclude that the spiral nature of the interior equilibrium implies cycles in crime and law enforcement. These occur because the concave shape of the inspector's cost function encourages 'extreme' levels of law enforcement. That is,  the inspector hardly enforces at all when crime is low as the average cost of enforcement is high. Conversely, enforcement is relatively strict when crime is high as the average cost of enforcement is relatively low. As a result, the inspector successively under and over-compensate. Meanwhile, the instability of these cycles stems from the fact that inspector discounts the future, therefore the optimal control allows crime to get out of hand in the long run.
\end{remark}

\section{Inspection with continuous strategy spaces}

In the previous sections, we assume that individuals can only make a binary choice between Violate and Comply. However, in many applications  individuals may be able to choose the extent of their criminal activity. For example, Kolokotsov and Malafeyev \cite{KM2010} analysed the case of a one-shot game between a tax authority and a single individual who can choose the amount of tax to conceal from the authority.

In this section we will extend  analysis in \cite{KM2010}  on inspection of a single individual by deriving the entire class of mixed strategy equilibria in this game. Then we will proceed to show that this class of equilibria also applies to games of population inspection.

\subsection{Inspection of a single individual}

\begin{figure}[ht]
\centering
\resizebox{1\textwidth}{!}{\includegraphics{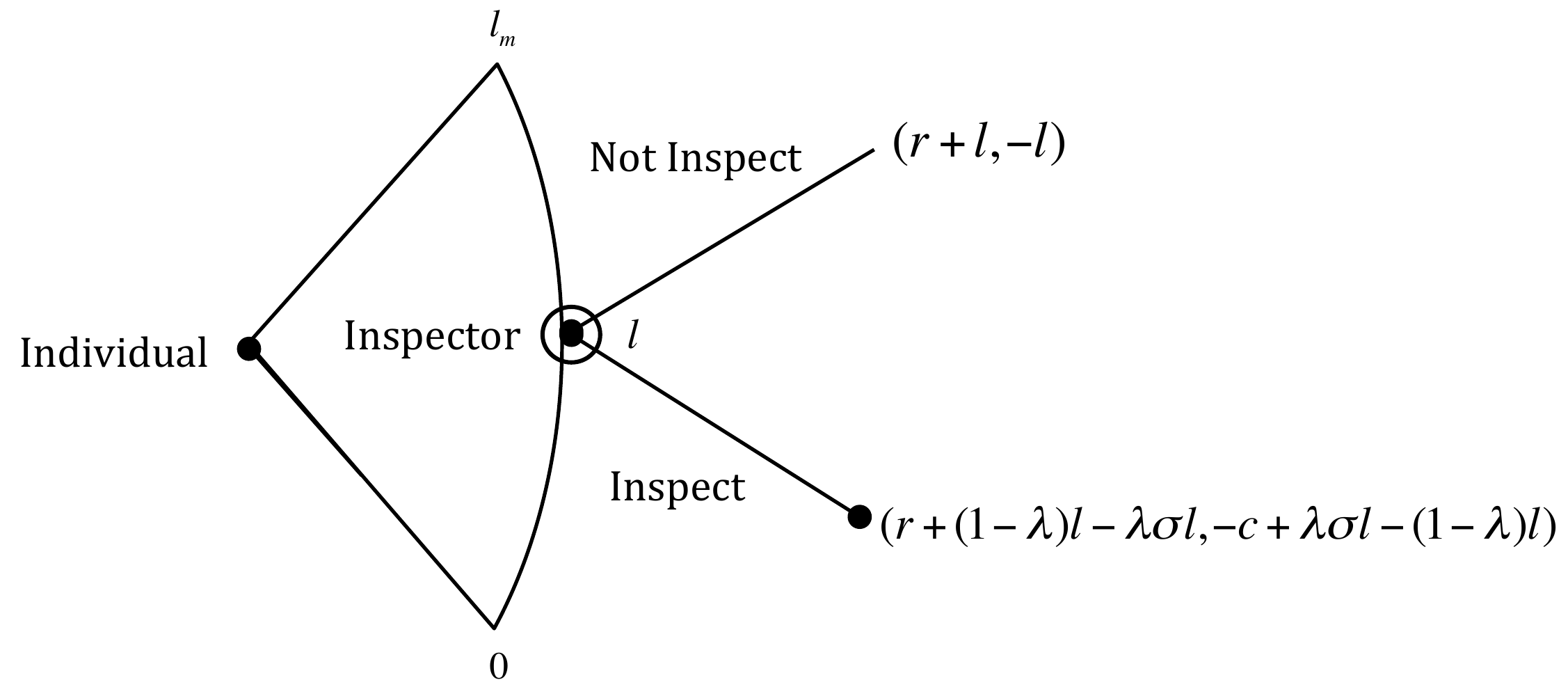} }
\caption{A single-individual inspection game with a continuous strategy space}
\label{A single-individual inspection game with continuous strategy choices}
\end{figure}

An inspection game of a single individual with a continuous strategy space has similar structure to the canonical game in a general form specified in Section \ref{Canonical games} (Table \ref{A general inspection game}), except that:
\begin{enumerate}
\item
the individual may now choose the extent of his criminal activity from the continuous strategy space $l\in[0,l_m]$ for a given $l_m>0$;\\[-1em]
\item
we assume that the fine levied by the inspector is proportional to the severity of the individual's crime, i.e. $f=\sigma l$ with $\sigma \in (0,1)$.
\end{enumerate}
It is worth stressing that, in this subsection we keep the assumption that the inspector chooses from the binary strategy space \{Inspect, Not Inspect\}. The payoffs to the inspector and the individual in the inspection game with continuous strategy space is summarised by the extensive form game in Figure \ref{A single-individual inspection game with continuous strategy choices}.

In order to ensure that this game does not degenerate into dominant strategies, we impose the following regulatory conditions on the payoffs: 
{\it 
\begin{equation*}
\hspace{-8em} \text {Assumption {\bf (A4)}} \quad\quad l_m>\frac{c}{\lambda(\sigma+1)} \,\text{ and } \,\,\lambda>\frac{1}{\sigma+1}.
\end{equation*}
}

The first inequality in assumption {\bf (A4)} states that the inspector would prefer to play Inspect when the individual chooses the maximum level of crime. The second inequality in assumption {\bf (A4)} states that  the individual would not wish to choose to commit any level of crime when the inspector is inspecting (otherwise maximum crime is always worthwhile).

Under assumption {\bf (A4)}, it is clear that there are no pure strategy Nash equilibria in this game. We shall search for mixed strategy equilibria where at least one player randomises over his pure strategies. Let $q\in[0,1]$ be the probability with which the inspector plays Inspect. We have the following result:
\begin{prop}\label{indiv  cont}
Let assumption {\bf (A4)} hold. In the inspection game with a continuous strategy space $[0,l_m]$, where $l_m>0$ is given, the equilibrium state is that the individual may play any randomisation over his strategy space provided that the expected level of crime is $l^*$ and the inspector inspects with the probability $q^*$, with
$$l^*=\frac{c}{\lambda(\sigma+1)}\, \text{ and } \, q^*=\frac{1}{\lambda(\sigma+1)}.$$
\end{prop}

\proof
Suppose for any probability of inspection $q\in[0,1]$, the individual will only randomise over two arbitrary pure strategies, $l_1,l_2\in [0,l_m]$ providing that these two pure strategies yield the same expected payoff (otherwise he would want to assign full probability to the strategy yielding the higher payoff).
By using the payoff table (Figure \ref{A single-individual inspection game with continuous strategy choices}), we get
\begin{equation}\label{Q}
\begin{split}
&q[r+(1-\lambda)l_1-\lambda \sigma l_1]+(1-q)[r+l_1]\\
=&q[r+(1-\lambda)l_2-\lambda \sigma l_2]+(1-q)[r+l_2].
\end{split}
\end{equation}
Solving equation for \eqref{Q} gives us the equilibrium probability of inspection
\begin{equation}\label{sol. Q}
q^*=\frac{1}{\lambda(\sigma +1)},
\end{equation}
which is independent of $l_1$ and $l_2$. This means, when the inspector plays the mixed strategy $q^*$ given by \eqref{sol. Q}, the individual is in fact willing to randomise over all of his pure strategies $l\in [0,l_m]$.

Suppose the individual's randomisation is given by some arbitrary probability density function $\phi(\cdot)$. Then, in order for the inspector to randomise, the expected payoff from playing Inspect must be equal to the expected payoff from playing Not Inspect, that is
$$\int_0^{l_m}(-c+\lambda \sigma l -(1-\lambda)l)\phi(l)dl=\int_0^{l_m}-l\phi(l)dl$$
which gives
\begin{equation}\label{l*}
\int_0^{l_m}l\phi(l)dl=\frac{c}{\lambda(\sigma+1)}:=l^*.
\end{equation}
Therefore, as long as the inspector plays the mixed strategy $q^*$ defined in \eqref{sol. Q}, the individual is willing to play any randomisation over all of his pure strategies $l\in [0,l_m]$; Meanwhile, as long as the individual plays an randomisation which results in the expected level of crime $l^*$ defined in \eqref{l*}, the inspector is willing to play any randomisation. If both players conform to these randomisations, the players' strategies are then mutually consistent best responses.
\qed

Proposition \ref{indiv  cont} implies a most interesting class of mixed strategy equilibria where the only constraint on the individual's equilibrium behaviour is the expected level of crime. Specifically, the individual may play any mixed strategy provided that his average level of crime is  $\frac{c}{\lambda(\sigma+1)}$.

\subsection{Population inspection with continuous strategy spaces}

We now consider  the inspector is responsible for the inspection of a large population of individuals. Specifically, we assume that there are $N$ individuals, each of whom chooses a level of crime $l$ from the continuous strategy space $[0,l_m]$. Let $\Phi(\cdot)$ be the probability density function of the distribution of crime in the population. The inspector chooses the level of law enforcement from a continuous strategy space, i.e. $q\in [0,1]$. The cost function of law enforcement for the inspector $F(\cdot)$ satisfies the assumption {\bf (A1)} in subsection \ref{Evolutionary inspection games}.

Then, given a rate of enforcement $q\in[0,1]$, the expected payoff for an individual from choosing the level of crime $l\in[0,l_m]$ is
$$\E[U_I(l,q)]=r+(1-\lambda q)l-\lambda q\sigma l.$$
Meanwhile, the expected payoff for the inspector is given by:
\begin{equation}\label{payoff au}
\E[U_A (q,\Phi(\cdot))]=-F(q)+\lambda qN\int_0^{l_m}\sigma l \Phi(l)dl-(1-\lambda q)N\int_0^{l_m}l\Phi(l)dl.
\end{equation}

\begin{prop}
 In an population inspection game, if the inspector chooses the level of law enforcement from a continuous strategy space, i.e. $q\in [0,1]$ and  each individual chooses a level of crime from a continuous strategy space, i.e. $l\in[0,l_m]$,  then any distribution of crime may occur in equilibrium if the expected level of crime is $L^*$ and the equilibrium probability of inspection is $q^*$, where
 $$L^*=\frac{F'(q^*)}{N\lambda (\sigma+1)}\, \text{ and } \,  q^*=\frac{1}{\lambda (\sigma+1)}.$$
\end{prop}

\proof
Following the same arguments in the proof of Proposition \ref{indiv  cont}, one can get $q^*$.  This quantity $q^*$ corresponds to the level of law enforcement at which individuals are willing to randomise over all pure strategies in the interval $[0,l_m]$.

Suppose the resulting crime probability density function is $\Phi(\cdot)$. The best response for the inspector  can be found by maximising \eqref{payoff au} over $q\in[0,1]$. The relevant first order condition is:
\begin{equation}\label{1order}
F'(q)=N \lambda (\sigma+1)\int_0^{l_m}l\Phi(l)dl.
\end{equation}
Setting $q=q^*$ in \eqref{1order} yields
\begin{equation}\label{Average L}
\int_0^{l_m}l\Phi(l)dl=\frac{F'(q^*)}{N\lambda (\sigma+1)}:=L^*.
\end{equation}
Therefore, $q^*$ is a best-response to $\Phi(\cdot)$ provided that the average level of crime is $L^*$ in \eqref{Average L}. Meanwhile, any distribution $\Phi(\cdot)$ with the average level of crime $L^*$ in \eqref{Average L} is a best-response to the level of law enforcement is $q^*$. Therefore $(L^*, q^*)$ constitutes a mixed strategy Nash equilibrium.\qed

This class of equilibria is particularly interesting for population-level inspection games, as it implies that populations with similar levels of law-enforcement can exhibit markedly different patterns of crime. In some populations, crime may be the preserve of a minority of relatively serious offenders, whereas in others, moderate levels of crime may be widespread. In terms of law enforcement, only the average level of crime matters.

\end{document}

%% file: rgb.tex
\definecolor{AliceBlue}{rgb}{0.94,0.97,1.00}
\definecolor{AntiqueWhite1}{rgb}{1.00,0.94,0.86}
\definecolor{AntiqueWhite2}{rgb}{0.93,0.87,0.80}
\definecolor{AntiqueWhite3}{rgb}{0.80,0.75,0.69}
\definecolor{AntiqueWhite4}{rgb}{0.55,0.51,0.47}
\definecolor{AntiqueWhite}{rgb}{0.98,0.92,0.84}
\definecolor{BlanchedAlmond}{rgb}{1.00,0.92,0.80}
\definecolor{BlueViolet}{rgb}{0.54,0.17,0.89}
\definecolor{CadetBlue1}{rgb}{0.60,0.96,1.00}
\definecolor{CadetBlue2}{rgb}{0.56,0.90,0.93}
\definecolor{CadetBlue3}{rgb}{0.48,0.77,0.80}
\definecolor{CadetBlue4}{rgb}{0.33,0.53,0.55}
\definecolor{CadetBlue}{rgb}{0.37,0.62,0.63}
\definecolor{CornflowerBlue}{rgb}{0.39,0.58,0.93}
\definecolor{DarkBlue}{rgb}{0.00,0.00,0.55}
\definecolor{DarkCyan}{rgb}{0.00,0.55,0.55}
\definecolor{DarkGoldenrod1}{rgb}{1.00,0.73,0.06}
\definecolor{DarkGoldenrod2}{rgb}{0.93,0.68,0.05}
\definecolor{DarkGoldenrod3}{rgb}{0.80,0.58,0.05}
\definecolor{DarkGoldenrod4}{rgb}{0.55,0.40,0.03}
\definecolor{DarkGoldenrod}{rgb}{0.72,0.53,0.04}
\definecolor{DarkGray}{rgb}{0.66,0.66,0.66}
\definecolor{DarkGreen}{rgb}{0.00,0.39,0.00}
\definecolor{DarkGrey}{rgb}{0.66,0.66,0.66}
\definecolor{DarkKhaki}{rgb}{0.74,0.72,0.42}
\definecolor{DarkMagenta}{rgb}{0.55,0.00,0.55}
\definecolor{DarkOliveGreen1}{rgb}{0.79,1.00,0.44}
\definecolor{DarkOliveGreen2}{rgb}{0.74,0.93,0.41}
\definecolor{DarkOliveGreen3}{rgb}{0.64,0.80,0.35}
\definecolor{DarkOliveGreen4}{rgb}{0.43,0.55,0.24}
\definecolor{DarkOliveGreen}{rgb}{0.33,0.42,0.18}
\definecolor{DarkOrange1}{rgb}{1.00,0.50,0.00}
\definecolor{DarkOrange2}{rgb}{0.93,0.46,0.00}
\definecolor{DarkOrange3}{rgb}{0.80,0.40,0.00}
\definecolor{DarkOrange4}{rgb}{0.55,0.27,0.00}
\definecolor{DarkOrange}{rgb}{1.00,0.55,0.00}
\definecolor{DarkOrchid1}{rgb}{0.75,0.24,1.00}
\definecolor{DarkOrchid2}{rgb}{0.70,0.23,0.93}
\definecolor{DarkOrchid3}{rgb}{0.60,0.20,0.80}
\definecolor{DarkOrchid4}{rgb}{0.41,0.13,0.55}
\definecolor{DarkOrchid}{rgb}{0.60,0.20,0.80}
\definecolor{DarkRed}{rgb}{0.55,0.00,0.00}
\definecolor{DarkSalmon}{rgb}{0.91,0.59,0.48}
\definecolor{DarkSeaGreen1}{rgb}{0.76,1.00,0.76}
\definecolor{DarkSeaGreen2}{rgb}{0.71,0.93,0.71}
\definecolor{DarkSeaGreen3}{rgb}{0.61,0.80,0.61}
\definecolor{DarkSeaGreen4}{rgb}{0.41,0.55,0.41}
\definecolor{DarkSeaGreen}{rgb}{0.56,0.74,0.56}
\definecolor{DarkSlateBlue}{rgb}{0.28,0.24,0.55}
\definecolor{DarkSlateGray1}{rgb}{0.59,1.00,1.00}
\definecolor{DarkSlateGray2}{rgb}{0.55,0.93,0.93}
\definecolor{DarkSlateGray3}{rgb}{0.47,0.80,0.80}
\definecolor{DarkSlateGray4}{rgb}{0.32,0.55,0.55}
\definecolor{DarkSlateGray}{rgb}{0.18,0.31,0.31}
\definecolor{DarkSlateGrey}{rgb}{0.18,0.31,0.31}
\definecolor{DarkTurquoise}{rgb}{0.00,0.81,0.82}
\definecolor{DarkViolet}{rgb}{0.58,0.00,0.83}
\definecolor{DeepPink1}{rgb}{1.00,0.08,0.58}
\definecolor{DeepPink2}{rgb}{0.93,0.07,0.54}
\definecolor{DeepPink3}{rgb}{0.80,0.06,0.46}
\definecolor{DeepPink4}{rgb}{0.55,0.04,0.31}
\definecolor{DeepPink}{rgb}{1.00,0.08,0.58}
\definecolor{DeepSkyBlue1}{rgb}{0.00,0.75,1.00}
\definecolor{DeepSkyBlue2}{rgb}{0.00,0.70,0.93}
\definecolor{DeepSkyBlue3}{rgb}{0.00,0.60,0.80}
\definecolor{DeepSkyBlue4}{rgb}{0.00,0.41,0.55}
\definecolor{DeepSkyBlue}{rgb}{0.00,0.75,1.00}
\definecolor{DimGray}{rgb}{0.41,0.41,0.41}
\definecolor{DimGrey}{rgb}{0.41,0.41,0.41}
\definecolor{DodgerBlue1}{rgb}{0.12,0.56,1.00}
\definecolor{DodgerBlue2}{rgb}{0.11,0.53,0.93}
\definecolor{DodgerBlue3}{rgb}{0.09,0.45,0.80}
\definecolor{DodgerBlue4}{rgb}{0.06,0.31,0.55}
\definecolor{DodgerBlue}{rgb}{0.12,0.56,1.00}
\definecolor{FloralWhite}{rgb}{1.00,0.98,0.94}
\definecolor{ForestGreen}{rgb}{0.13,0.55,0.13}
\definecolor{GhostWhite}{rgb}{0.97,0.97,1.00}
\definecolor{GreenYellow}{rgb}{0.68,1.00,0.18}
\definecolor{HotPink1}{rgb}{1.00,0.43,0.71}
\definecolor{HotPink2}{rgb}{0.93,0.42,0.65}
\definecolor{HotPink3}{rgb}{0.80,0.38,0.56}
\definecolor{HotPink4}{rgb}{0.55,0.23,0.38}
\definecolor{HotPink}{rgb}{1.00,0.41,0.71}
\definecolor{IndianRed1}{rgb}{1.00,0.42,0.42}
\definecolor{IndianRed2}{rgb}{0.93,0.39,0.39}
\definecolor{IndianRed3}{rgb}{0.80,0.33,0.33}
\definecolor{IndianRed4}{rgb}{0.55,0.23,0.23}
\definecolor{IndianRed}{rgb}{0.80,0.36,0.36}
\definecolor{LavenderBlush1}{rgb}{1.00,0.94,0.96}
\definecolor{LavenderBlush2}{rgb}{0.93,0.88,0.90}
\definecolor{LavenderBlush3}{rgb}{0.80,0.76,0.77}
\definecolor{LavenderBlush4}{rgb}{0.55,0.51,0.53}
\definecolor{LavenderBlush}{rgb}{1.00,0.94,0.96}
\definecolor{LawnGreen}{rgb}{0.49,0.99,0.00}
\definecolor{LemonChiffon1}{rgb}{1.00,0.98,0.80}
\definecolor{LemonChiffon2}{rgb}{0.93,0.91,0.75}
\definecolor{LemonChiffon3}{rgb}{0.80,0.79,0.65}
\definecolor{LemonChiffon4}{rgb}{0.55,0.54,0.44}
\definecolor{LemonChiffon}{rgb}{1.00,0.98,0.80}
\definecolor{LightBlue1}{rgb}{0.75,0.94,1.00}
\definecolor{LightBlue2}{rgb}{0.70,0.87,0.93}
\definecolor{LightBlue3}{rgb}{0.60,0.75,0.80}
\definecolor{LightBlue4}{rgb}{0.41,0.51,0.55}
\definecolor{LightBlue}{rgb}{0.68,0.85,0.90}
\definecolor{LightCoral}{rgb}{0.94,0.50,0.50}
\definecolor{LightCyan1}{rgb}{0.88,1.00,1.00}
\definecolor{LightCyan2}{rgb}{0.82,0.93,0.93}
\definecolor{LightCyan3}{rgb}{0.71,0.80,0.80}
\definecolor{LightCyan4}{rgb}{0.48,0.55,0.55}
\definecolor{LightCyan}{rgb}{0.88,1.00,1.00}
\definecolor{LightGoldenrod1}{rgb}{1.00,0.93,0.55}
\definecolor{LightGoldenrod2}{rgb}{0.93,0.86,0.51}
\definecolor{LightGoldenrod3}{rgb}{0.80,0.75,0.44}
\definecolor{LightGoldenrod4}{rgb}{0.55,0.51,0.30}
\definecolor{LightGoldenrodYellow}{rgb}{0.98,0.98,0.82}
\definecolor{LightGoldenrod}{rgb}{0.93,0.87,0.51}
\definecolor{LightGray}{rgb}{0.83,0.83,0.83}
\definecolor{LightGreen}{rgb}{0.56,0.93,0.56}
\definecolor{LightGrey}{rgb}{0.83,0.83,0.83}
\definecolor{LightPink1}{rgb}{1.00,0.68,0.73}
\definecolor{LightPink2}{rgb}{0.93,0.64,0.68}
\definecolor{LightPink3}{rgb}{0.80,0.55,0.58}
\definecolor{LightPink4}{rgb}{0.55,0.37,0.40}
\definecolor{LightPink}{rgb}{1.00,0.71,0.76}
\definecolor{LightSalmon1}{rgb}{1.00,0.63,0.48}
\definecolor{LightSalmon2}{rgb}{0.93,0.58,0.45}
\definecolor{LightSalmon3}{rgb}{0.80,0.51,0.38}
\definecolor{LightSalmon4}{rgb}{0.55,0.34,0.26}
\definecolor{LightSalmon}{rgb}{1.00,0.63,0.48}
\definecolor{LightSeaGreen}{rgb}{0.13,0.70,0.67}
\definecolor{LightSkyBlue1}{rgb}{0.69,0.89,1.00}
\definecolor{LightSkyBlue2}{rgb}{0.64,0.83,0.93}
\definecolor{LightSkyBlue3}{rgb}{0.55,0.71,0.80}
\definecolor{LightSkyBlue4}{rgb}{0.38,0.48,0.55}
\definecolor{LightSkyBlue}{rgb}{0.53,0.81,0.98}
\definecolor{LightSlateBlue}{rgb}{0.52,0.44,1.00}
\definecolor{LightSlateGray}{rgb}{0.47,0.53,0.60}
\definecolor{LightSlateGrey}{rgb}{0.47,0.53,0.60}
\definecolor{LightSteelBlue1}{rgb}{0.79,0.88,1.00}
\definecolor{LightSteelBlue2}{rgb}{0.74,0.82,0.93}
\definecolor{LightSteelBlue3}{rgb}{0.64,0.71,0.80}
\definecolor{LightSteelBlue4}{rgb}{0.43,0.48,0.55}
\definecolor{LightSteelBlue}{rgb}{0.69,0.77,0.87}
\definecolor{LightYellow1}{rgb}{1.00,1.00,0.88}
\definecolor{LightYellow2}{rgb}{0.93,0.93,0.82}
\definecolor{LightYellow3}{rgb}{0.80,0.80,0.71}
\definecolor{LightYellow4}{rgb}{0.55,0.55,0.48}
\definecolor{LightYellow}{rgb}{1.00,1.00,0.88}
\definecolor{LimeGreen}{rgb}{0.20,0.80,0.20}
\definecolor{MediumAquamarine}{rgb}{0.40,0.80,0.67}
\definecolor{MediumBlue}{rgb}{0.00,0.00,0.80}
\definecolor{MediumOrchid1}{rgb}{0.88,0.40,1.00}
\definecolor{MediumOrchid2}{rgb}{0.82,0.37,0.93}
\definecolor{MediumOrchid3}{rgb}{0.71,0.32,0.80}
\definecolor{MediumOrchid4}{rgb}{0.48,0.22,0.55}
\definecolor{MediumOrchid}{rgb}{0.73,0.33,0.83}
\definecolor{MediumPurple1}{rgb}{0.67,0.51,1.00}
\definecolor{MediumPurple2}{rgb}{0.62,0.47,0.93}
\definecolor{MediumPurple3}{rgb}{0.54,0.41,0.80}
\definecolor{MediumPurple4}{rgb}{0.36,0.28,0.55}
\definecolor{MediumPurple}{rgb}{0.58,0.44,0.86}
\definecolor{MediumSeaGreen}{rgb}{0.24,0.70,0.44}
\definecolor{MediumSlateBlue}{rgb}{0.48,0.41,0.93}
\definecolor{MediumSpringGreen}{rgb}{0.00,0.98,0.60}
\definecolor{MediumTurquoise}{rgb}{0.28,0.82,0.80}
\definecolor{MediumVioletRed}{rgb}{0.78,0.08,0.52}
\definecolor{MidnightBlue}{rgb}{0.10,0.10,0.44}
\definecolor{MintCream}{rgb}{0.96,1.00,0.98}
\definecolor{MistyRose1}{rgb}{1.00,0.89,0.88}
\definecolor{MistyRose2}{rgb}{0.93,0.84,0.82}
\definecolor{MistyRose3}{rgb}{0.80,0.72,0.71}
\definecolor{MistyRose4}{rgb}{0.55,0.49,0.48}
\definecolor{MistyRose}{rgb}{1.00,0.89,0.88}
\definecolor{NavajoWhite1}{rgb}{1.00,0.87,0.68}
\definecolor{NavajoWhite2}{rgb}{0.93,0.81,0.63}
\definecolor{NavajoWhite3}{rgb}{0.80,0.70,0.55}
\definecolor{NavajoWhite4}{rgb}{0.55,0.47,0.37}
\definecolor{NavajoWhite}{rgb}{1.00,0.87,0.68}
\definecolor{NavyBlue}{rgb}{0.00,0.00,0.50}
\definecolor{OldLace}{rgb}{0.99,0.96,0.90}
\definecolor{OliveDrab1}{rgb}{0.75,1.00,0.24}
\definecolor{OliveDrab2}{rgb}{0.70,0.93,0.23}
\definecolor{OliveDrab3}{rgb}{0.60,0.80,0.20}
\definecolor{OliveDrab4}{rgb}{0.41,0.55,0.13}
\definecolor{OliveDrab}{rgb}{0.42,0.56,0.14}
\definecolor{OrangeRed1}{rgb}{1.00,0.27,0.00}
\definecolor{OrangeRed2}{rgb}{0.93,0.25,0.00}
\definecolor{OrangeRed3}{rgb}{0.80,0.22,0.00}
\definecolor{OrangeRed4}{rgb}{0.55,0.15,0.00}
\definecolor{OrangeRed}{rgb}{1.00,0.27,0.00}
\definecolor{PaleGoldenrod}{rgb}{0.93,0.91,0.67}
\definecolor{PaleGreen1}{rgb}{0.60,1.00,0.60}
\definecolor{PaleGreen2}{rgb}{0.56,0.93,0.56}
\definecolor{PaleGreen3}{rgb}{0.49,0.80,0.49}
\definecolor{PaleGreen4}{rgb}{0.33,0.55,0.33}
\definecolor{PaleGreen}{rgb}{0.60,0.98,0.60}
\definecolor{PaleTurquoise1}{rgb}{0.73,1.00,1.00}
\definecolor{PaleTurquoise2}{rgb}{0.68,0.93,0.93}
\definecolor{PaleTurquoise3}{rgb}{0.59,0.80,0.80}
\definecolor{PaleTurquoise4}{rgb}{0.40,0.55,0.55}
\definecolor{PaleTurquoise}{rgb}{0.69,0.93,0.93}
\definecolor{PaleVioletRed1}{rgb}{1.00,0.51,0.67}
\definecolor{PaleVioletRed2}{rgb}{0.93,0.47,0.62}
\definecolor{PaleVioletRed3}{rgb}{0.80,0.41,0.54}
\definecolor{PaleVioletRed4}{rgb}{0.55,0.28,0.36}
\definecolor{PaleVioletRed}{rgb}{0.86,0.44,0.58}
\definecolor{PapayaWhip}{rgb}{1.00,0.94,0.84}
\definecolor{PeachPuff1}{rgb}{1.00,0.85,0.73}
\definecolor{PeachPuff2}{rgb}{0.93,0.80,0.68}
\definecolor{PeachPuff3}{rgb}{0.80,0.69,0.58}
\definecolor{PeachPuff4}{rgb}{0.55,0.47,0.40}
\definecolor{PeachPuff}{rgb}{1.00,0.85,0.73}
\definecolor{PowderBlue}{rgb}{0.69,0.88,0.90}
\definecolor{RosyBrown1}{rgb}{1.00,0.76,0.76}
\definecolor{RosyBrown2}{rgb}{0.93,0.71,0.71}
\definecolor{RosyBrown3}{rgb}{0.80,0.61,0.61}
\definecolor{RosyBrown4}{rgb}{0.55,0.41,0.41}
\definecolor{RosyBrown}{rgb}{0.74,0.56,0.56}
\definecolor{RoyalBlue1}{rgb}{0.28,0.46,1.00}
\definecolor{RoyalBlue2}{rgb}{0.26,0.43,0.93}
\definecolor{RoyalBlue3}{rgb}{0.23,0.37,0.80}
\definecolor{RoyalBlue4}{rgb}{0.15,0.25,0.55}
\definecolor{RoyalBlue}{rgb}{0.25,0.41,0.88}
\definecolor{SaddleBrown}{rgb}{0.55,0.27,0.07}
\definecolor{SandyBrown}{rgb}{0.96,0.64,0.38}
\definecolor{SeaGreen1}{rgb}{0.33,1.00,0.62}
\definecolor{SeaGreen2}{rgb}{0.31,0.93,0.58}
\definecolor{SeaGreen3}{rgb}{0.26,0.80,0.50}
\definecolor{SeaGreen4}{rgb}{0.18,0.55,0.34}
\definecolor{SeaGreen}{rgb}{0.18,0.55,0.34}
\definecolor{SkyBlue1}{rgb}{0.53,0.81,1.00}
\definecolor{SkyBlue2}{rgb}{0.49,0.75,0.93}
\definecolor{SkyBlue3}{rgb}{0.42,0.65,0.80}
\definecolor{SkyBlue4}{rgb}{0.29,0.44,0.55}
\definecolor{SkyBlue}{rgb}{0.53,0.81,0.92}
\definecolor{SlateBlue1}{rgb}{0.51,0.44,1.00}
\definecolor{SlateBlue2}{rgb}{0.48,0.40,0.93}
\definecolor{SlateBlue3}{rgb}{0.41,0.35,0.80}
\definecolor{SlateBlue4}{rgb}{0.28,0.24,0.55}
\definecolor{SlateBlue}{rgb}{0.42,0.35,0.80}
\definecolor{SlateGray1}{rgb}{0.78,0.89,1.00}
\definecolor{SlateGray2}{rgb}{0.73,0.83,0.93}
\definecolor{SlateGray3}{rgb}{0.62,0.71,0.80}
\definecolor{SlateGray4}{rgb}{0.42,0.48,0.55}
\definecolor{SlateGray}{rgb}{0.44,0.50,0.56}
\definecolor{SlateGrey}{rgb}{0.44,0.50,0.56}
\definecolor{SpringGreen1}{rgb}{0.00,1.00,0.50}
\definecolor{SpringGreen2}{rgb}{0.00,0.93,0.46}
\definecolor{SpringGreen3}{rgb}{0.00,0.80,0.40}
\definecolor{SpringGreen4}{rgb}{0.00,0.55,0.27}
\definecolor{SpringGreen}{rgb}{0.00,1.00,0.50}
\definecolor{SteelBlue1}{rgb}{0.39,0.72,1.00}
\definecolor{SteelBlue2}{rgb}{0.36,0.67,0.93}
\definecolor{SteelBlue3}{rgb}{0.31,0.58,0.80}
\definecolor{SteelBlue4}{rgb}{0.21,0.39,0.55}
\definecolor{SteelBlue}{rgb}{0.27,0.51,0.71}
\definecolor{VioletRed1}{rgb}{1.00,0.24,0.59}
\definecolor{VioletRed2}{rgb}{0.93,0.23,0.55}
\definecolor{VioletRed3}{rgb}{0.80,0.20,0.47}
\definecolor{VioletRed4}{rgb}{0.55,0.13,0.32}
\definecolor{VioletRed}{rgb}{0.82,0.13,0.56}
\definecolor{WhiteSmoke}{rgb}{0.96,0.96,0.96}
\definecolor{YellowGreen}{rgb}{0.60,0.80,0.20}
\definecolor{aliceblue}{rgb}{0.94,0.97,1.00}
\definecolor{antiquewhite}{rgb}{0.98,0.92,0.84}
\definecolor{aquamarine1}{rgb}{0.50,1.00,0.83}
\definecolor{aquamarine2}{rgb}{0.46,0.93,0.78}
\definecolor{aquamarine3}{rgb}{0.40,0.80,0.67}
\definecolor{aquamarine4}{rgb}{0.27,0.55,0.45}
\definecolor{aquamarine}{rgb}{0.50,1.00,0.83}
\definecolor{azure1}{rgb}{0.94,1.00,1.00}
\definecolor{azure2}{rgb}{0.88,0.93,0.93}
\definecolor{azure3}{rgb}{0.76,0.80,0.80}
\definecolor{azure4}{rgb}{0.51,0.55,0.55}
\definecolor{azure}{rgb}{0.94,1.00,1.00}
\definecolor{beige}{rgb}{0.96,0.96,0.86}
\definecolor{bisque1}{rgb}{1.00,0.89,0.77}
\definecolor{bisque2}{rgb}{0.93,0.84,0.72}
\definecolor{bisque3}{rgb}{0.80,0.72,0.62}
\definecolor{bisque4}{rgb}{0.55,0.49,0.42}
\definecolor{bisque}{rgb}{1.00,0.89,0.77}
\definecolor{black}{rgb}{0.00,0.00,0.00}
\definecolor{blanchedalmond}{rgb}{1.00,0.92,0.80}
\definecolor{blue1}{rgb}{0.00,0.00,1.00}
\definecolor{blue2}{rgb}{0.00,0.00,0.93}
\definecolor{blue3}{rgb}{0.00,0.00,0.80}
\definecolor{blue4}{rgb}{0.00,0.00,0.55}
\definecolor{blueviolet}{rgb}{0.54,0.17,0.89}
\definecolor{blue}{rgb}{0.00,0.00,1.00}
\definecolor{brown1}{rgb}{1.00,0.25,0.25}
\definecolor{brown2}{rgb}{0.93,0.23,0.23}
\definecolor{brown3}{rgb}{0.80,0.20,0.20}
\definecolor{brown4}{rgb}{0.55,0.14,0.14}
\definecolor{brown}{rgb}{0.65,0.16,0.16}
\definecolor{burlywood1}{rgb}{1.00,0.83,0.61}
\definecolor{burlywood2}{rgb}{0.93,0.77,0.57}
\definecolor{burlywood3}{rgb}{0.80,0.67,0.49}
\definecolor{burlywood4}{rgb}{0.55,0.45,0.33}
\definecolor{burlywood}{rgb}{0.87,0.72,0.53}
\definecolor{cadetblue}{rgb}{0.37,0.62,0.63}
\definecolor{chartreuse1}{rgb}{0.50,1.00,0.00}
\definecolor{chartreuse2}{rgb}{0.46,0.93,0.00}
\definecolor{chartreuse3}{rgb}{0.40,0.80,0.00}
\definecolor{chartreuse4}{rgb}{0.27,0.55,0.00}
\definecolor{chartreuse}{rgb}{0.50,1.00,0.00}
\definecolor{chocolate1}{rgb}{1.00,0.50,0.14}
\definecolor{chocolate2}{rgb}{0.93,0.46,0.13}
\definecolor{chocolate3}{rgb}{0.80,0.40,0.11}
\definecolor{chocolate4}{rgb}{0.55,0.27,0.07}
\definecolor{chocolate}{rgb}{0.82,0.41,0.12}
\definecolor{coral1}{rgb}{1.00,0.45,0.34}
\definecolor{coral2}{rgb}{0.93,0.42,0.31}
\definecolor{coral3}{rgb}{0.80,0.36,0.27}
\definecolor{coral4}{rgb}{0.55,0.24,0.18}
\definecolor{coral}{rgb}{1.00,0.50,0.31}
\definecolor{cornflowerblue}{rgb}{0.39,0.58,0.93}
\definecolor{cornsilk1}{rgb}{1.00,0.97,0.86}
\definecolor{cornsilk2}{rgb}{0.93,0.91,0.80}
\definecolor{cornsilk3}{rgb}{0.80,0.78,0.69}
\definecolor{cornsilk4}{rgb}{0.55,0.53,0.47}
\definecolor{cornsilk}{rgb}{1.00,0.97,0.86}
\definecolor{cyan1}{rgb}{0.00,1.00,1.00}
\definecolor{cyan2}{rgb}{0.00,0.93,0.93}
\definecolor{cyan3}{rgb}{0.00,0.80,0.80}
\definecolor{cyan4}{rgb}{0.00,0.55,0.55}
\definecolor{cyan}{rgb}{0.00,1.00,1.00}
\definecolor{darkblue}{rgb}{0.00,0.00,0.55}
\definecolor{darkcyan}{rgb}{0.00,0.55,0.55}
\definecolor{darkgoldenrod}{rgb}{0.72,0.53,0.04}
\definecolor{darkgray}{rgb}{0.66,0.66,0.66}
\definecolor{darkgreen}{rgb}{0.00,0.39,0.00}
\definecolor{darkgrey}{rgb}{0.66,0.66,0.66}
\definecolor{darkkhaki}{rgb}{0.74,0.72,0.42}
\definecolor{darkmagenta}{rgb}{0.55,0.00,0.55}
\definecolor{darkolive}{rgb}{0.33,0.42,0.18}
\definecolor{darkorange}{rgb}{1.00,0.55,0.00}
\definecolor{darkorchid}{rgb}{0.60,0.20,0.80}
\definecolor{darkred}{rgb}{0.55,0.00,0.00}
\definecolor{darksalmon}{rgb}{0.91,0.59,0.48}
\definecolor{darksea}{rgb}{0.56,0.74,0.56}
\definecolor{darkslate}{rgb}{0.18,0.31,0.31}
\definecolor{darkslate}{rgb}{0.18,0.31,0.31}
\definecolor{darkslate}{rgb}{0.28,0.24,0.55}
\definecolor{darkturquoise}{rgb}{0.00,0.81,0.82}
\definecolor{darkviolet}{rgb}{0.58,0.00,0.83}
\definecolor{deeppink}{rgb}{1.00,0.08,0.58}
\definecolor{deepsky}{rgb}{0.00,0.75,1.00}
\definecolor{dimgray}{rgb}{0.41,0.41,0.41}
\definecolor{dimgrey}{rgb}{0.41,0.41,0.41}
\definecolor{dodgerblue}{rgb}{0.12,0.56,1.00}
\definecolor{firebrick1}{rgb}{1.00,0.19,0.19}
\definecolor{firebrick2}{rgb}{0.93,0.17,0.17}
\definecolor{firebrick3}{rgb}{0.80,0.15,0.15}
\definecolor{firebrick4}{rgb}{0.55,0.10,0.10}
\definecolor{firebrick}{rgb}{0.70,0.13,0.13}
\definecolor{floralwhite}{rgb}{1.00,0.98,0.94}
\definecolor{forestgreen}{rgb}{0.13,0.55,0.13}
\definecolor{gainsboro}{rgb}{0.86,0.86,0.86}
\definecolor{ghostwhite}{rgb}{0.97,0.97,1.00}
\definecolor{gold1}{rgb}{1.00,0.84,0.00}
\definecolor{gold2}{rgb}{0.93,0.79,0.00}
\definecolor{gold3}{rgb}{0.80,0.68,0.00}
\definecolor{gold4}{rgb}{0.55,0.46,0.00}
\definecolor{goldenrod1}{rgb}{1.00,0.76,0.15}
\definecolor{goldenrod2}{rgb}{0.93,0.71,0.13}
\definecolor{goldenrod3}{rgb}{0.80,0.61,0.11}
\definecolor{goldenrod4}{rgb}{0.55,0.41,0.08}
\definecolor{goldenrod}{rgb}{0.85,0.65,0.13}
\definecolor{gold}{rgb}{1.00,0.84,0.00}
\definecolor{gray0}{rgb}{0.00,0.00,0.00}
\definecolor{gray100}{rgb}{1.00,1.00,1.00}
\definecolor{gray10}{rgb}{0.10,0.10,0.10}
\definecolor{gray11}{rgb}{0.11,0.11,0.11}
\definecolor{gray12}{rgb}{0.12,0.12,0.12}
\definecolor{gray13}{rgb}{0.13,0.13,0.13}
\definecolor{gray14}{rgb}{0.14,0.14,0.14}
\definecolor{gray15}{rgb}{0.15,0.15,0.15}
\definecolor{gray16}{rgb}{0.16,0.16,0.16}
\definecolor{gray17}{rgb}{0.17,0.17,0.17}
\definecolor{gray18}{rgb}{0.18,0.18,0.18}
\definecolor{gray19}{rgb}{0.19,0.19,0.19}
\definecolor{gray1}{rgb}{0.01,0.01,0.01}
\definecolor{gray20}{rgb}{0.20,0.20,0.20}
\definecolor{gray21}{rgb}{0.21,0.21,0.21}
\definecolor{gray22}{rgb}{0.22,0.22,0.22}
\definecolor{gray23}{rgb}{0.23,0.23,0.23}
\definecolor{gray24}{rgb}{0.24,0.24,0.24}
\definecolor{gray25}{rgb}{0.25,0.25,0.25}
\definecolor{gray26}{rgb}{0.26,0.26,0.26}
\definecolor{gray27}{rgb}{0.27,0.27,0.27}
\definecolor{gray28}{rgb}{0.28,0.28,0.28}
\definecolor{gray29}{rgb}{0.29,0.29,0.29}
\definecolor{gray2}{rgb}{0.02,0.02,0.02}
\definecolor{gray30}{rgb}{0.30,0.30,0.30}
\definecolor{gray31}{rgb}{0.31,0.31,0.31}
\definecolor{gray32}{rgb}{0.32,0.32,0.32}
\definecolor{gray33}{rgb}{0.33,0.33,0.33}
\definecolor{gray34}{rgb}{0.34,0.34,0.34}
\definecolor{gray35}{rgb}{0.35,0.35,0.35}
\definecolor{gray36}{rgb}{0.36,0.36,0.36}
\definecolor{gray37}{rgb}{0.37,0.37,0.37}
\definecolor{gray38}{rgb}{0.38,0.38,0.38}
\definecolor{gray39}{rgb}{0.39,0.39,0.39}
\definecolor{gray3}{rgb}{0.03,0.03,0.03}
\definecolor{gray40}{rgb}{0.40,0.40,0.40}
\definecolor{gray41}{rgb}{0.41,0.41,0.41}
\definecolor{gray42}{rgb}{0.42,0.42,0.42}
\definecolor{gray43}{rgb}{0.43,0.43,0.43}
\definecolor{gray44}{rgb}{0.44,0.44,0.44}
\definecolor{gray45}{rgb}{0.45,0.45,0.45}
\definecolor{gray46}{rgb}{0.46,0.46,0.46}
\definecolor{gray47}{rgb}{0.47,0.47,0.47}
\definecolor{gray48}{rgb}{0.48,0.48,0.48}
\definecolor{gray49}{rgb}{0.49,0.49,0.49}
\definecolor{gray4}{rgb}{0.04,0.04,0.04}
\definecolor{gray50}{rgb}{0.50,0.50,0.50}
\definecolor{gray51}{rgb}{0.51,0.51,0.51}
\definecolor{gray52}{rgb}{0.52,0.52,0.52}
\definecolor{gray53}{rgb}{0.53,0.53,0.53}
\definecolor{gray54}{rgb}{0.54,0.54,0.54}
\definecolor{gray55}{rgb}{0.55,0.55,0.55}
\definecolor{gray56}{rgb}{0.56,0.56,0.56}
\definecolor{gray57}{rgb}{0.57,0.57,0.57}
\definecolor{gray58}{rgb}{0.58,0.58,0.58}
\definecolor{gray59}{rgb}{0.59,0.59,0.59}
\definecolor{gray5}{rgb}{0.05,0.05,0.05}
\definecolor{gray60}{rgb}{0.60,0.60,0.60}
\definecolor{gray61}{rgb}{0.61,0.61,0.61}
\definecolor{gray62}{rgb}{0.62,0.62,0.62}
\definecolor{gray63}{rgb}{0.63,0.63,0.63}
\definecolor{gray64}{rgb}{0.64,0.64,0.64}
\definecolor{gray65}{rgb}{0.65,0.65,0.65}
\definecolor{gray66}{rgb}{0.66,0.66,0.66}
\definecolor{gray67}{rgb}{0.67,0.67,0.67}
\definecolor{gray68}{rgb}{0.68,0.68,0.68}
\definecolor{gray69}{rgb}{0.69,0.69,0.69}
\definecolor{gray6}{rgb}{0.06,0.06,0.06}
\definecolor{gray70}{rgb}{0.70,0.70,0.70}
\definecolor{gray71}{rgb}{0.71,0.71,0.71}
\definecolor{gray72}{rgb}{0.72,0.72,0.72}
\definecolor{gray73}{rgb}{0.73,0.73,0.73}
\definecolor{gray74}{rgb}{0.74,0.74,0.74}
\definecolor{gray75}{rgb}{0.75,0.75,0.75}
\definecolor{gray76}{rgb}{0.76,0.76,0.76}
\definecolor{gray77}{rgb}{0.77,0.77,0.77}
\definecolor{gray78}{rgb}{0.78,0.78,0.78}
\definecolor{gray79}{rgb}{0.79,0.79,0.79}
\definecolor{gray7}{rgb}{0.07,0.07,0.07}
\definecolor{gray80}{rgb}{0.80,0.80,0.80}
\definecolor{gray81}{rgb}{0.81,0.81,0.81}
\definecolor{gray82}{rgb}{0.82,0.82,0.82}
\definecolor{gray83}{rgb}{0.83,0.83,0.83}
\definecolor{gray84}{rgb}{0.84,0.84,0.84}
\definecolor{gray85}{rgb}{0.85,0.85,0.85}
\definecolor{gray86}{rgb}{0.86,0.86,0.86}
\definecolor{gray87}{rgb}{0.87,0.87,0.87}
\definecolor{gray88}{rgb}{0.88,0.88,0.88}
\definecolor{gray89}{rgb}{0.89,0.89,0.89}
\definecolor{gray8}{rgb}{0.08,0.08,0.08}
\definecolor{gray90}{rgb}{0.90,0.90,0.90}
\definecolor{gray91}{rgb}{0.91,0.91,0.91}
\definecolor{gray92}{rgb}{0.92,0.92,0.92}
\definecolor{gray93}{rgb}{0.93,0.93,0.93}
\definecolor{gray94}{rgb}{0.94,0.94,0.94}
\definecolor{gray95}{rgb}{0.95,0.95,0.95}
\definecolor{gray96}{rgb}{0.96,0.96,0.96}
\definecolor{gray97}{rgb}{0.97,0.97,0.97}
\definecolor{gray98}{rgb}{0.98,0.98,0.98}
\definecolor{gray99}{rgb}{0.99,0.99,0.99}
\definecolor{gray9}{rgb}{0.09,0.09,0.09}
\definecolor{gray}{rgb}{0.75,0.75,0.75}
\definecolor{green1}{rgb}{0.00,1.00,0.00}
\definecolor{green2}{rgb}{0.00,0.93,0.00}
\definecolor{green3}{rgb}{0.00,0.80,0.00}
\definecolor{green4}{rgb}{0.00,0.55,0.00}
\definecolor{greenyellow}{rgb}{0.68,1.00,0.18}
\definecolor{green}{rgb}{0.00,1.00,0.00}
\definecolor{grey0}{rgb}{0.00,0.00,0.00}
\definecolor{grey100}{rgb}{1.00,1.00,1.00}
\definecolor{grey10}{rgb}{0.10,0.10,0.10}
\definecolor{grey11}{rgb}{0.11,0.11,0.11}
\definecolor{grey12}{rgb}{0.12,0.12,0.12}
\definecolor{grey13}{rgb}{0.13,0.13,0.13}
\definecolor{grey14}{rgb}{0.14,0.14,0.14}
\definecolor{grey15}{rgb}{0.15,0.15,0.15}
\definecolor{grey16}{rgb}{0.16,0.16,0.16}
\definecolor{grey17}{rgb}{0.17,0.17,0.17}
\definecolor{grey18}{rgb}{0.18,0.18,0.18}
\definecolor{grey19}{rgb}{0.19,0.19,0.19}
\definecolor{grey1}{rgb}{0.01,0.01,0.01}
\definecolor{grey20}{rgb}{0.20,0.20,0.20}
\definecolor{grey21}{rgb}{0.21,0.21,0.21}
\definecolor{grey22}{rgb}{0.22,0.22,0.22}
\definecolor{grey23}{rgb}{0.23,0.23,0.23}
\definecolor{grey24}{rgb}{0.24,0.24,0.24}
\definecolor{grey25}{rgb}{0.25,0.25,0.25}
\definecolor{grey26}{rgb}{0.26,0.26,0.26}
\definecolor{grey27}{rgb}{0.27,0.27,0.27}
\definecolor{grey28}{rgb}{0.28,0.28,0.28}
\definecolor{grey29}{rgb}{0.29,0.29,0.29}
\definecolor{grey2}{rgb}{0.02,0.02,0.02}
\definecolor{grey30}{rgb}{0.30,0.30,0.30}
\definecolor{grey31}{rgb}{0.31,0.31,0.31}
\definecolor{grey32}{rgb}{0.32,0.32,0.32}
\definecolor{grey33}{rgb}{0.33,0.33,0.33}
\definecolor{grey34}{rgb}{0.34,0.34,0.34}
\definecolor{grey35}{rgb}{0.35,0.35,0.35}
\definecolor{grey36}{rgb}{0.36,0.36,0.36}
\definecolor{grey37}{rgb}{0.37,0.37,0.37}
\definecolor{grey38}{rgb}{0.38,0.38,0.38}
\definecolor{grey39}{rgb}{0.39,0.39,0.39}
\definecolor{grey3}{rgb}{0.03,0.03,0.03}
\definecolor{grey40}{rgb}{0.40,0.40,0.40}
\definecolor{grey41}{rgb}{0.41,0.41,0.41}
\definecolor{grey42}{rgb}{0.42,0.42,0.42}
\definecolor{grey43}{rgb}{0.43,0.43,0.43}
\definecolor{grey44}{rgb}{0.44,0.44,0.44}
\definecolor{grey45}{rgb}{0.45,0.45,0.45}
\definecolor{grey46}{rgb}{0.46,0.46,0.46}
\definecolor{grey47}{rgb}{0.47,0.47,0.47}
\definecolor{grey48}{rgb}{0.48,0.48,0.48}
\definecolor{grey49}{rgb}{0.49,0.49,0.49}
\definecolor{grey4}{rgb}{0.04,0.04,0.04}
\definecolor{grey50}{rgb}{0.50,0.50,0.50}
\definecolor{grey51}{rgb}{0.51,0.51,0.51}
\definecolor{grey52}{rgb}{0.52,0.52,0.52}
\definecolor{grey53}{rgb}{0.53,0.53,0.53}
\definecolor{grey54}{rgb}{0.54,0.54,0.54}
\definecolor{grey55}{rgb}{0.55,0.55,0.55}
\definecolor{grey56}{rgb}{0.56,0.56,0.56}
\definecolor{grey57}{rgb}{0.57,0.57,0.57}
\definecolor{grey58}{rgb}{0.58,0.58,0.58}
\definecolor{grey59}{rgb}{0.59,0.59,0.59}
\definecolor{grey5}{rgb}{0.05,0.05,0.05}
\definecolor{grey60}{rgb}{0.60,0.60,0.60}
\definecolor{grey61}{rgb}{0.61,0.61,0.61}
\definecolor{grey62}{rgb}{0.62,0.62,0.62}
\definecolor{grey63}{rgb}{0.63,0.63,0.63}
\definecolor{grey64}{rgb}{0.64,0.64,0.64}
\definecolor{grey65}{rgb}{0.65,0.65,0.65}
\definecolor{grey66}{rgb}{0.66,0.66,0.66}
\definecolor{grey67}{rgb}{0.67,0.67,0.67}
\definecolor{grey68}{rgb}{0.68,0.68,0.68}
\definecolor{grey69}{rgb}{0.69,0.69,0.69}
\definecolor{grey6}{rgb}{0.06,0.06,0.06}
\definecolor{grey70}{rgb}{0.70,0.70,0.70}
\definecolor{grey71}{rgb}{0.71,0.71,0.71}
\definecolor{grey72}{rgb}{0.72,0.72,0.72}
\definecolor{grey73}{rgb}{0.73,0.73,0.73}
\definecolor{grey74}{rgb}{0.74,0.74,0.74}
\definecolor{grey75}{rgb}{0.75,0.75,0.75}
\definecolor{grey76}{rgb}{0.76,0.76,0.76}
\definecolor{grey77}{rgb}{0.77,0.77,0.77}
\definecolor{grey78}{rgb}{0.78,0.78,0.78}
\definecolor{grey79}{rgb}{0.79,0.79,0.79}
\definecolor{grey7}{rgb}{0.07,0.07,0.07}
\definecolor{grey80}{rgb}{0.80,0.80,0.80}
\definecolor{grey81}{rgb}{0.81,0.81,0.81}
\definecolor{grey82}{rgb}{0.82,0.82,0.82}
\definecolor{grey83}{rgb}{0.83,0.83,0.83}
\definecolor{grey84}{rgb}{0.84,0.84,0.84}
\definecolor{grey85}{rgb}{0.85,0.85,0.85}
\definecolor{grey86}{rgb}{0.86,0.86,0.86}
\definecolor{grey87}{rgb}{0.87,0.87,0.87}
\definecolor{grey88}{rgb}{0.88,0.88,0.88}
\definecolor{grey89}{rgb}{0.89,0.89,0.89}
\definecolor{grey8}{rgb}{0.08,0.08,0.08}
\definecolor{grey90}{rgb}{0.90,0.90,0.90}
\definecolor{grey91}{rgb}{0.91,0.91,0.91}
\definecolor{grey92}{rgb}{0.92,0.92,0.92}
\definecolor{grey93}{rgb}{0.93,0.93,0.93}
\definecolor{grey94}{rgb}{0.94,0.94,0.94}
\definecolor{grey95}{rgb}{0.95,0.95,0.95}
\definecolor{grey96}{rgb}{0.96,0.96,0.96}
\definecolor{grey97}{rgb}{0.97,0.97,0.97}
\definecolor{grey98}{rgb}{0.98,0.98,0.98}
\definecolor{grey99}{rgb}{0.99,0.99,0.99}
\definecolor{grey9}{rgb}{0.09,0.09,0.09}
\definecolor{grey}{rgb}{0.75,0.75,0.75}
\definecolor{honeydew1}{rgb}{0.94,1.00,0.94}
\definecolor{honeydew2}{rgb}{0.88,0.93,0.88}
\definecolor{honeydew3}{rgb}{0.76,0.80,0.76}
\definecolor{honeydew4}{rgb}{0.51,0.55,0.51}
\definecolor{honeydew}{rgb}{0.94,1.00,0.94}
\definecolor{hotpink}{rgb}{1.00,0.41,0.71}
\definecolor{indianred}{rgb}{0.80,0.36,0.36}
\definecolor{ivory1}{rgb}{1.00,1.00,0.94}
\definecolor{ivory2}{rgb}{0.93,0.93,0.88}
\definecolor{ivory3}{rgb}{0.80,0.80,0.76}
\definecolor{ivory4}{rgb}{0.55,0.55,0.51}
\definecolor{ivory}{rgb}{1.00,1.00,0.94}
\definecolor{khaki1}{rgb}{1.00,0.96,0.56}
\definecolor{khaki2}{rgb}{0.93,0.90,0.52}
\definecolor{khaki3}{rgb}{0.80,0.78,0.45}
\definecolor{khaki4}{rgb}{0.55,0.53,0.31}
\definecolor{khaki}{rgb}{0.94,0.90,0.55}
\definecolor{lavenderblush}{rgb}{1.00,0.94,0.96}
\definecolor{lavender}{rgb}{0.90,0.90,0.98}
\definecolor{lawngreen}{rgb}{0.49,0.99,0.00}
\definecolor{lemonchiffon}{rgb}{1.00,0.98,0.80}
\definecolor{lightblue}{rgb}{0.68,0.85,0.90}
\definecolor{lightcoral}{rgb}{0.94,0.50,0.50}
\definecolor{lightcyan}{rgb}{0.88,1.00,1.00}
\definecolor{lightgoldenrod}{rgb}{0.93,0.87,0.51}
\definecolor{lightgoldenrod}{rgb}{0.98,0.98,0.82}
\definecolor{lightgray}{rgb}{0.83,0.83,0.83}
\definecolor{lightgreen}{rgb}{0.56,0.93,0.56}
\definecolor{lightgrey}{rgb}{0.83,0.83,0.83}
\definecolor{lightpink}{rgb}{1.00,0.71,0.76}
\definecolor{lightsalmon}{rgb}{1.00,0.63,0.48}
\definecolor{lightsea}{rgb}{0.13,0.70,0.67}
\definecolor{lightsky}{rgb}{0.53,0.81,0.98}
\definecolor{lightslate}{rgb}{0.47,0.53,0.60}
\definecolor{lightslate}{rgb}{0.47,0.53,0.60}
\definecolor{lightslate}{rgb}{0.52,0.44,1.00}
\definecolor{lightsteel}{rgb}{0.69,0.77,0.87}
\definecolor{lightyellow}{rgb}{1.00,1.00,0.88}
\definecolor{limegreen}{rgb}{0.20,0.80,0.20}
\definecolor{linen}{rgb}{0.98,0.94,0.90}
\definecolor{magenta1}{rgb}{1.00,0.00,1.00}
\definecolor{magenta2}{rgb}{0.93,0.00,0.93}
\definecolor{magenta3}{rgb}{0.80,0.00,0.80}
\definecolor{magenta4}{rgb}{0.55,0.00,0.55}
\definecolor{magenta}{rgb}{1.00,0.00,1.00}
\definecolor{maroon1}{rgb}{1.00,0.20,0.70}
\definecolor{maroon2}{rgb}{0.93,0.19,0.65}
\definecolor{maroon3}{rgb}{0.80,0.16,0.56}
\definecolor{maroon4}{rgb}{0.55,0.11,0.38}
\definecolor{maroon}{rgb}{0.69,0.19,0.38}
\definecolor{mediumaquamarine}{rgb}{0.40,0.80,0.67}
\definecolor{mediumblue}{rgb}{0.00,0.00,0.80}
\definecolor{mediumorchid}{rgb}{0.73,0.33,0.83}
\definecolor{mediumpurple}{rgb}{0.58,0.44,0.86}
\definecolor{mediumsea}{rgb}{0.24,0.70,0.44}
\definecolor{mediumslate}{rgb}{0.48,0.41,0.93}
\definecolor{mediumspring}{rgb}{0.00,0.98,0.60}
\definecolor{mediumturquoise}{rgb}{0.28,0.82,0.80}
\definecolor{mediumviolet}{rgb}{0.78,0.08,0.52}
\definecolor{midnightblue}{rgb}{0.10,0.10,0.44}
\definecolor{mintcream}{rgb}{0.96,1.00,0.98}
\definecolor{mistyrose}{rgb}{1.00,0.89,0.88}
\definecolor{moccasin}{rgb}{1.00,0.89,0.71}
\definecolor{navajowhite}{rgb}{1.00,0.87,0.68}
\definecolor{navyblue}{rgb}{0.00,0.00,0.50}
\definecolor{navy}{rgb}{0.00,0.00,0.50}
\definecolor{oldlace}{rgb}{0.99,0.96,0.90}
\definecolor{olivedrab}{rgb}{0.42,0.56,0.14}
\definecolor{orange1}{rgb}{1.00,0.65,0.00}
\definecolor{orange2}{rgb}{0.93,0.60,0.00}
\definecolor{orange3}{rgb}{0.80,0.52,0.00}
\definecolor{orange4}{rgb}{0.55,0.35,0.00}
\definecolor{orangered}{rgb}{1.00,0.27,0.00}
\definecolor{orange}{rgb}{1.00,0.65,0.00}
\definecolor{orchid1}{rgb}{1.00,0.51,0.98}
\definecolor{orchid2}{rgb}{0.93,0.48,0.91}
\definecolor{orchid3}{rgb}{0.80,0.41,0.79}
\definecolor{orchid4}{rgb}{0.55,0.28,0.54}
\definecolor{orchid}{rgb}{0.85,0.44,0.84}
\definecolor{palegoldenrod}{rgb}{0.93,0.91,0.67}
\definecolor{palegreen}{rgb}{0.60,0.98,0.60}
\definecolor{paleturquoise}{rgb}{0.69,0.93,0.93}
\definecolor{paleviolet}{rgb}{0.86,0.44,0.58}
\definecolor{papayawhip}{rgb}{1.00,0.94,0.84}
\definecolor{peachpuff}{rgb}{1.00,0.85,0.73}
\definecolor{peru}{rgb}{0.80,0.52,0.25}
\definecolor{pink1}{rgb}{1.00,0.71,0.77}
\definecolor{pink2}{rgb}{0.93,0.66,0.72}
\definecolor{pink3}{rgb}{0.80,0.57,0.62}
\definecolor{pink4}{rgb}{0.55,0.39,0.42}
\definecolor{pink}{rgb}{1.00,0.75,0.80}
\definecolor{plum1}{rgb}{1.00,0.73,1.00}
\definecolor{plum2}{rgb}{0.93,0.68,0.93}
\definecolor{plum3}{rgb}{0.80,0.59,0.80}
\definecolor{plum4}{rgb}{0.55,0.40,0.55}
\definecolor{plum}{rgb}{0.87,0.63,0.87}
\definecolor{powderblue}{rgb}{0.69,0.88,0.90}
\definecolor{purple1}{rgb}{0.61,0.19,1.00}
\definecolor{purple2}{rgb}{0.57,0.17,0.93}
\definecolor{purple3}{rgb}{0.49,0.15,0.80}
\definecolor{purple4}{rgb}{0.33,0.10,0.55}
\definecolor{purple}{rgb}{0.63,0.13,0.94}
\definecolor{red1}{rgb}{1.00,0.00,0.00}
\definecolor{red2}{rgb}{0.93,0.00,0.00}
\definecolor{red3}{rgb}{0.80,0.00,0.00}
\definecolor{red4}{rgb}{0.55,0.00,0.00}
\definecolor{red}{rgb}{1.00,0.00,0.00}
\definecolor{rosybrown}{rgb}{0.74,0.56,0.56}
\definecolor{royalblue}{rgb}{0.25,0.41,0.88}
\definecolor{saddlebrown}{rgb}{0.55,0.27,0.07}
\definecolor{salmon1}{rgb}{1.00,0.55,0.41}
\definecolor{salmon2}{rgb}{0.93,0.51,0.38}
\definecolor{salmon3}{rgb}{0.80,0.44,0.33}
\definecolor{salmon4}{rgb}{0.55,0.30,0.22}
\definecolor{salmon}{rgb}{0.98,0.50,0.45}
\definecolor{sandybrown}{rgb}{0.96,0.64,0.38}
\definecolor{seagreen}{rgb}{0.18,0.55,0.34}
\definecolor{seashell1}{rgb}{1.00,0.96,0.93}
\definecolor{seashell2}{rgb}{0.93,0.90,0.87}
\definecolor{seashell3}{rgb}{0.80,0.77,0.75}
\definecolor{seashell4}{rgb}{0.55,0.53,0.51}
\definecolor{seashell}{rgb}{1.00,0.96,0.93}
\definecolor{sienna1}{rgb}{1.00,0.51,0.28}
\definecolor{sienna2}{rgb}{0.93,0.47,0.26}
\definecolor{sienna3}{rgb}{0.80,0.41,0.22}
\definecolor{sienna4}{rgb}{0.55,0.28,0.15}
\definecolor{sienna}{rgb}{0.63,0.32,0.18}
\definecolor{skyblue}{rgb}{0.53,0.81,0.92}
\definecolor{slateblue}{rgb}{0.42,0.35,0.80}
\definecolor{slategray}{rgb}{0.44,0.50,0.56}
\definecolor{slategrey}{rgb}{0.44,0.50,0.56}
\definecolor{snow1}{rgb}{1.00,0.98,0.98}
\definecolor{snow2}{rgb}{0.93,0.91,0.91}
\definecolor{snow3}{rgb}{0.80,0.79,0.79}
\definecolor{snow4}{rgb}{0.55,0.54,0.54}
\definecolor{snow}{rgb}{1.00,0.98,0.98}
\definecolor{springgreen}{rgb}{0.00,1.00,0.50}
\definecolor{steelblue}{rgb}{0.27,0.51,0.71}
\definecolor{tan1}{rgb}{1.00,0.65,0.31}
\definecolor{tan2}{rgb}{0.93,0.60,0.29}
\definecolor{tan3}{rgb}{0.80,0.52,0.25}
\definecolor{tan4}{rgb}{0.55,0.35,0.17}
\definecolor{tan}{rgb}{0.82,0.71,0.55}
\definecolor{thistle1}{rgb}{1.00,0.88,1.00}
\definecolor{thistle2}{rgb}{0.93,0.82,0.93}
\definecolor{thistle3}{rgb}{0.80,0.71,0.80}
\definecolor{thistle4}{rgb}{0.55,0.48,0.55}
\definecolor{thistle}{rgb}{0.85,0.75,0.85}
\definecolor{tomato1}{rgb}{1.00,0.39,0.28}
\definecolor{tomato2}{rgb}{0.93,0.36,0.26}
\definecolor{tomato3}{rgb}{0.80,0.31,0.22}
\definecolor{tomato4}{rgb}{0.55,0.21,0.15}
\definecolor{tomato}{rgb}{1.00,0.39,0.28}
\definecolor{turquoise1}{rgb}{0.00,0.96,1.00}
\definecolor{turquoise2}{rgb}{0.00,0.90,0.93}
\definecolor{turquoise3}{rgb}{0.00,0.77,0.80}
\definecolor{turquoise4}{rgb}{0.00,0.53,0.55}
\definecolor{turquoise}{rgb}{0.25,0.88,0.82}
\definecolor{violetred}{rgb}{0.82,0.13,0.56}
\definecolor{violet}{rgb}{0.93,0.51,0.93}
\definecolor{wheat1}{rgb}{1.00,0.91,0.73}
\definecolor{wheat2}{rgb}{0.93,0.85,0.68}
\definecolor{wheat3}{rgb}{0.80,0.73,0.59}
\definecolor{wheat4}{rgb}{0.55,0.49,0.40}
\definecolor{wheat}{rgb}{0.96,0.87,0.70}
\definecolor{whitesmoke}{rgb}{0.96,0.96,0.96}
\definecolor{white}{rgb}{1.00,1.00,1.00}
\definecolor{yellow1}{rgb}{1.00,1.00,0.00}
\definecolor{yellow2}{rgb}{0.93,0.93,0.00}
\definecolor{yellow3}{rgb}{0.80,0.80,0.00}
\definecolor{yellow4}{rgb}{0.55,0.55,0.00}
\definecolor{yellowgreen}{rgb}{0.60,0.80,0.20}
\definecolor{yellow}{rgb}{1.00,1.00,0.00}